\documentclass[onecolumn,11pt]{article}
\usepackage[top=1in, bottom=1in, left=1in, right=1in]{geometry}
\usepackage{amsmath,amsfonts,amscd,amssymb}
\usepackage{graphicx}
\usepackage{epstopdf}
\usepackage{overpic}
\usepackage{cancel}
\usepackage{rotating}
\usepackage{url}
\usepackage{caption}
\usepackage{color}
\usepackage{rotating}
\usepackage{multirow}
\usepackage{wrapfig}
\usepackage{mathtools}
\usepackage{subeqnarray}
\usepackage{setspace}
\usepackage{palatino} 
\usepackage[numbers,sort&compress]{natbib}
\usepackage{array}
\usepackage{tikz}
\PassOptionsToPackage{linktocpage}{hyperref}
\newcolumntype{+}{!{\vrule width 2pt}}

\usepackage{booktabs,longtable,tabu} 
\usepackage{multirow} 
\usepackage{float}
\usepackage{makecell}

\newlength\savedwidth

\usepackage{microtype}
\DisableLigatures[f]{encoding = *, family = * }
\usepackage[font=small,margin=0.25in,labelsep={colon},justification=justified,]{caption}
\usepackage{subcaption}

\usepackage[bottom,flushmargin,hang,multiple]{footmisc}
\usepackage{lipsum}
\newcommand\blfootnote[1]{%
  \begingroup
  \renewcommand\thefootnote{}\footnote{#1}%
  \addtocounter{footnote}{-1}%
  \endgroup
}

\definecolor{header1}{cmyk}{0,0,0,1}

\title{\LARGE{\vspace{-.55in}\textbf{Randomized methods to characterize large-scale vortical flow network}}}

\author{\normalsize{Zhe Bai$^{1*}$, N. Benjamin Erichson$^2$, Muralikrishnan Gopalakrishnan Meena$^3$},
\\ \normalsize{Kunihiko Taira$^3$, Steven L. Brunton$^{4*}$}\\
\footnotesize{$^1$ Computational Research Division, Lawrence Berkeley National Laboratory, Berkeley, CA 94720, United States}\\
\footnotesize{$^2$ Department of Applied Mathematics, University of Washington, Seattle, WA 98195, United States}\\
\footnotesize{$^3$ Department of Mechanical and Aerospace Engineering, University of California, Los Angeles, CA 90095, United States}\\
\footnotesize{$^4$ Department of Mechanical Engineering, University of Washington, Seattle, WA 98195, United States}
\vspace{-.2in}
}
\date{}

\begin{document}
\maketitle
\blfootnote{$^*$ Corresponding author (zhebai@lbl.gov).}
\vspace{-.2in}

\begin{abstract}
We demonstrate the effective use of randomized methods for linear algebra to perform network-based analysis of complex vortical flows. Network theoretic approaches can reveal the connectivity structures among a set of vortical elements and analyze their collective dynamics. These approaches have recently been generalized to analyze high-dimensional turbulent flows, for which network computations can become prohibitively expensive.  In this work, we propose efficient methods to approximate network quantities, such as the leading eigendecomposition of the adjacency matrix, using randomized methods. 
Specifically, we use the Nystr\"om method to approximate the leading eigenvalues and eigenvectors, achieving significant computational savings and reduced memory requirements.
The effectiveness of the proposed technique is demonstrated on two high-dimensional flow fields: two-dimensional flow past an airfoil and two-dimensional turbulence.
We find that quasi-uniform column sampling outperforms uniform column sampling, while both feature the same computational complexity.

\vspace{.15in}

\end{abstract}

\section{Introduction}
Fluid dynamics is a rich and challenging field at the intersection of physics and engineering.  
There is still a tremendous amount that we do not understand about turbulence, yet working fluids are at the heart of nearly every major industry, including health, defense, transportation, and energy.  
One cannot overestimate the significance and impact that a better understanding of fluid flows would have in our ability to predict and manipulate their behavior in the real world.  
The challenge of modeling and controlling fluids stems from the fact that they are nonlinear with complex multi-scale interactions over a large range of spatial and temporal scales. 
Moreover, data from experiments and simulations are generally represented as exceedingly high-dimensional measurements that may obscure underlying patterns.   
With advances in simulation capabilities and measurement techniques, the volume and quality of such data are rapidly increasing.  

Despite the complexity of the governing equations and the overwhelming volume of data, it is often observed that flows evolve on a low-dimensional attractor defined by a few dominant coherent structures~\cite{HLBR_turb}. 
There are a wealth of modal decomposition techniques to mine and characterize these structures from experimental data and numerical solvers~\cite{Taira2017aiaa,taira2019modal}.  
The majority of techniques, such as proper orthogonal decomposition (POD)\cite{HLBR_turb,Noack2003jfm} and dynamic mode decomposition~\cite{Schmid2010jfm,Rowley2009jfm,Kutz2016book} are fundamentally linear, as are many of the standard techniques in control theory~\cite{dp:book}.  
Linear control has been widely applied for flow control~\cite{Bagheri2009amr,Brunton2015amr,Sipp2016amr}, for example to stabilize boundary layers.  
However, many flows are fundamentally nonlinear, limiting the broad use of linear control theory.  
Low-order nonlinear models may be obtained by Galerkin projection of the Navier Stokes equations onto POD modes, and these models have had great success for uncovering underling mechanisms that drive flows~\cite{Noack2003jfm}. 
Regardless, there are issues with Galerkin models, such as instability and mode deformation with changing parameters and boundary conditions~\cite{Carlberg2017jcp,Loiseau2017jfm,Loiseau2018jfm}, limiting the success of POD-Galerkin modeling for turbulence. 

Recently, network theoretic approaches have been increasingly leveraged to analyze complex, fluid flow systems~\cite{Nair2015JFM,murugesan2015combustion,Taira2016JFM,schlueter2017jfm,meena2018network,murayama2018characterization,iacobello2019lagrangian}. 
Network science~\cite{Newman10} characterizes the structure and dynamics on a graph consisting of nodes and the edges connecting them. 
It is possible to cast a fluid flow in this context, where each grid cell of vorticity is a node, and the induced velocity from the Biot-Savart interaction between each node establishes the edge connections~\cite{Nair2015JFM}.  
In this way, individual packets of vorticity interact and evolve according to physics that is encoded in the graph.  
Network analysis provides a complementary perspective to classical techniques in fluid dynamics, especially for flow control.  
There have been many powerful advances in network control theory~\cite{Newman10,Mesbahi2010book} surrounding multi-agent systems~\cite{Rahmani:SIAMJCO09,Mesbahi2010book} in the past two decades.  
Multi-agent control is designed to be local, efficient to operate, and scalable to extremely large systems, such as the internet~\cite{Low2002ieeecs,Doyle2005pnas} or the electric grid~\cite{Susuki2011jns}.  
The multi-agent control and underlying network dynamics may also be strongly nonlinear, and these controllers are built to handle time-delays and communication failures, which are typical limiting factors for robust performance in classical control~\cite{dp:book}.
In particular, networks are often characterized by a large collection of elements, represented by nodes on a graph, that each execute their own set of local protocols in response to external stimulus.  
This analogy holds quite well for a number of large networked dynamical systems, including animals flocking~\cite{Leonard2001cdc,Olfati2006ieeetac}, multi-robotic cooperative control systems~\cite{Balch1998ieeetac}, sensor networks~\cite{Cortes2002ieeera,Leonard2007pieee}, biological regulatory networks~\cite{Milo2002science,Luscombe2004nature}, and the internet~\cite{Low2002ieeecs,Doyle2005pnas}, to name a few.  
A long-term goal of characterizing vortical networks is the eventual application of multi-agent control to \emph{school} turbulence into a beneficial configuration for an engineering advantage.  
In contrast, past closed-loop flow control efforts have commonly involved applying linear control techniques to a suitable reduced-order model of the fluid~\cite{Brunton2015amr}.

\subsection{Motivation}
The application of network theory in fluid dynamics faces the challenge of an exceedingly large number of degrees of freedom (nodes) for a fully turbulent flow, leading to an even larger adjacency matrix of edges. 
Indeed, the adjacency matrix for a vortical flow network scales as the square of the number of fluid grid points, quickly becoming intractable to store and analyze using classical techniques from linear algebra. 
Recent approaches have been developed to manage this complexity, including building a graph on the \emph{modes} of a flow~\cite{Nair2018pre}, given by dominant fluid coherent structures from POD, rather than the \emph{nodes}, given by fluid grid elements.
However, there is still interest in developing scalable methods to directly analyze the vortical network in the original ambient measurement space, where nodes are grid elements. 
Randomized methods for linear algebra provide an emerging alternative to efficiently compute an approximate  eigendecomposition of large-scale matrices. 
These methods work with a reduced representation, a so-called \emph{sketch}, of the input matrix that captures the essential spectral information. 
This sketch is obtained either through \emph{random sampling} or by random \emph{projections}.
The \emph{sketch} can then be used to efficiently compute a desired low-rank approximation, such as the singular value decomposition. 
For a comprehensive survey and implementation details see~\cite{halko2011finding,DrineasRNLA,kannan2017randomized,woodruff2014sketching,erichson_jss}.

\subsection{Contributions}
\begin{figure}[t]
\hspace{-.25in}\begin{overpic}[trim = {0, 0, 0, 0}, clip, width=1.09\textwidth]{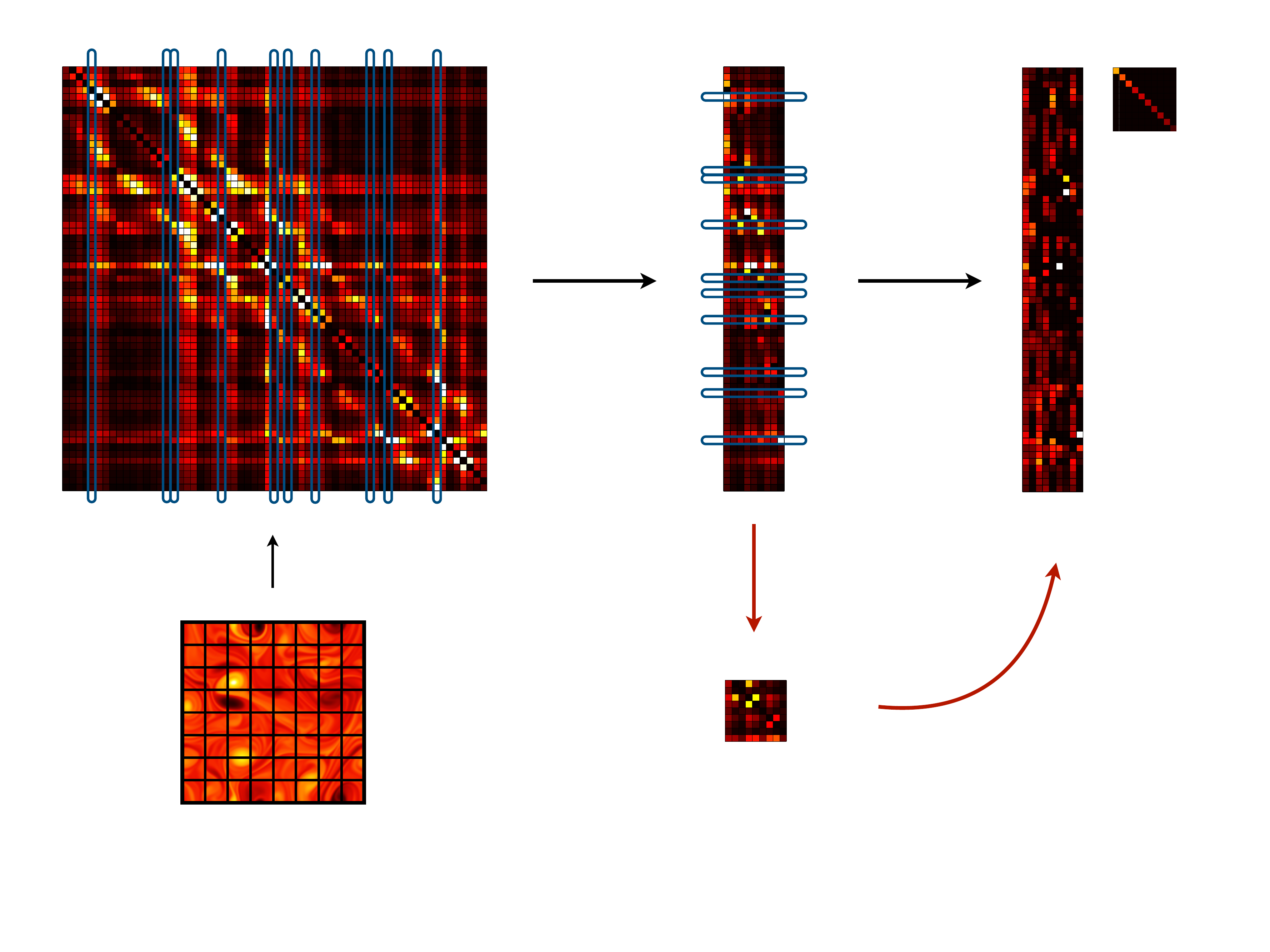}
		\put(14, 9){\textit{Input Flow Field}}
		\put(12.8, 73){{Adjacency Matrix}}
		\put(42.8, 57){\small{Column}}
		\put(42.0, 55){\small{Sampling}}
		
		\put(42.0, 29.0){\small{Row Sampling}}
		\put(56.6, 73){\small{Sketch}}
		\put(80, 73){\small{Eigenvectors}}
		\put(79, 70.6){\small{\& Eigenvalues}}
		\put(69, 55){\small{Path \bf{I}}}
		\put(64, 16.0){\small Path {\bf{II}} (Nystr\"om Method)}
	\end{overpic}
	\vspace{-.7in}
\caption{Identification of dominant eigenvalues and eigenvectors of the adjacency matrix obtained from a fluid flow field using randomized methods. Path $\mathbf{I}$ randomly samples columns of the adjacency matrix to approximate the leading eigenvalues and eigenvectors; path $\mathbf{II}$ uses both column and row sampling to form a sketch. Note, it is not required to explicitly construct the full adjacency matrix.}
\label{fig1}
\end{figure}
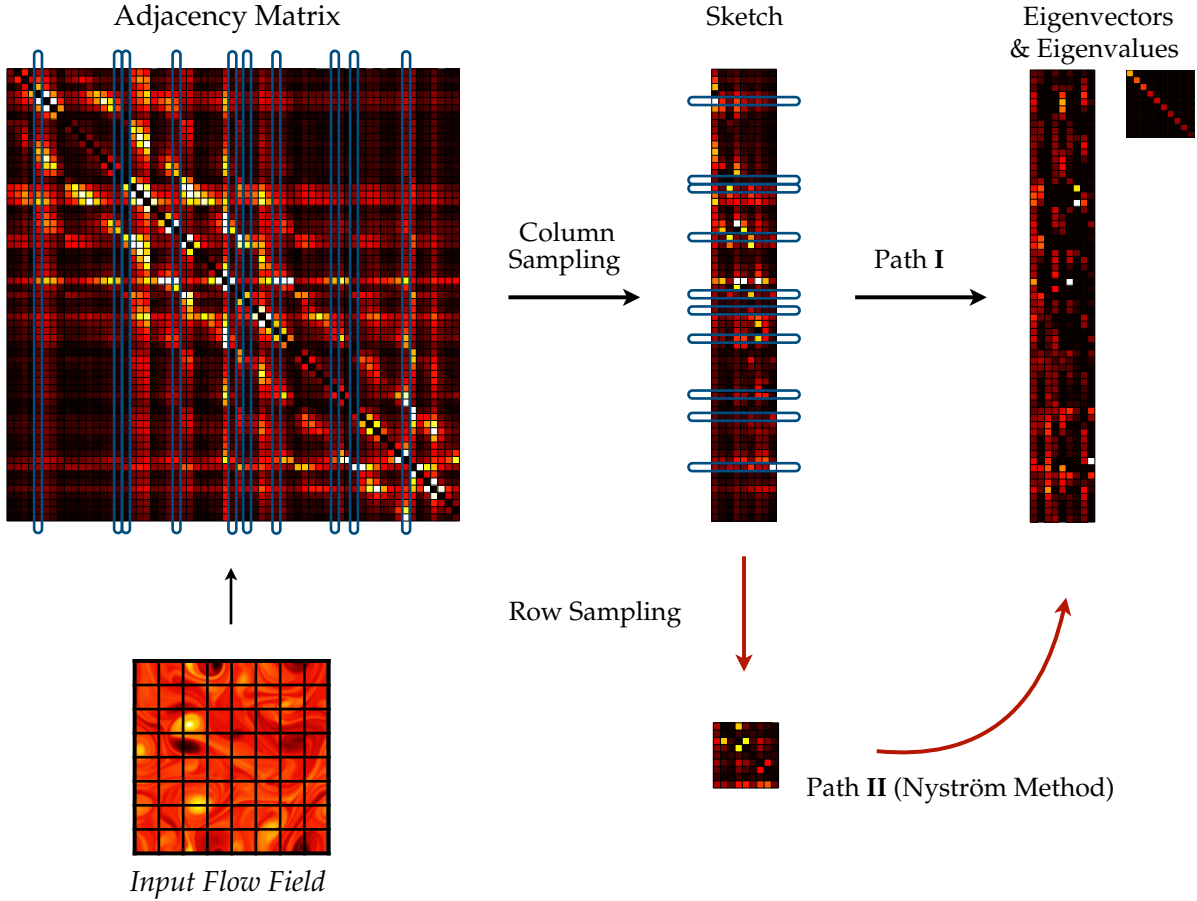

In this work, we use randomized methods to analyze large networks arising in fluid mechanics.  
Specifically we are interested in using randomized methods to obtain a qualitatively correct view of the dominant graph structures, which can then be used for the downstream tasks of sensor and actuator placement, optimization, and control.  
This application domain is the subject of the present work, where we explore the use of random sampling and randomized linear algebra to generate qualitatively accurate estimates of dominant flow structures for cases where the adjacency matrix is exceedingly large.  
We compare both uniform and quasi-uniform sampling, and demonstrate this approach on two high-dimensional vortical flow networks.  
A high-level principle sketch is provided in Figure~\ref{fig1}.

\section{Vortical interaction network}
Recently, network analysis has been employed to represent and analyze vortical interactions in flows~\cite{Nair2015JFM,Taira2016JFM,meena2018network,Nair2018pre,Nair2019jfm}. 
In Nair and Taira~\cite{Nair2015JFM}, spectral sparsification was used to identify a low-dimensional skeleton underlying the high-dimensional, nonlinear fluid dynamics.  
Taira et al.~\cite{Taira2016JFM} later demonstrated scale-free behavior of the vortex interaction network for two-dimensional isotropic turbulence.  
Network analysis has also been used by Meena et al.~\cite{meena2018network} for community detection and force prediction in an unsteady wake. 
More recently, interaction networks based on the energy transfer between \emph{modes} has been used for modeling and control~\cite{Nair2018pre,Nair2019jfm}.
These studies provide compelling evidence that the emerging network perspective may complement well-established flow modeling techniques.

{A graph (network) $\mathcal{G}$ consists of three components: a set of vertices (nodes) $\mathcal{V}$, a set of edges $\mathcal{E}$ that connect these vertices, and the associated edge weights $\mathcal{W}$. In the present work, the vortical elements in a flow field represent the nodes of the vortical flow network, following the formulation of Nair and Taira~\cite{Nair2015JFM}. The interaction among vortical elements can be quantified by the induced velocity.  Following the Biot-Savart law, the induced velocity from element $i$ to $j$ in a two-dimensional flow is given by
\begin{equation}
    u_{i\rightarrow j} = \frac{\Gamma_i}{2\pi d_{ij}}
\end{equation}
where $\Gamma_i = \omega_i \Delta x \Delta y$ is the strength of vortical element $i$ with vorticity $\omega_i$ and grid area $\Delta x \Delta y$, and $d_{ij}$ is the Euclidean distance between the vortical elements. 
This induced velocity represents the edge weights, and may be used to construct an adjacency matrix 
\begin{align}\label{e:adj}
\mathbf{A}_{ij} = 
\begin{cases}
	 \frac{1}{2} \left( |u_{i \rightarrow j}| + |u_{j \rightarrow i}|\right) & \text{if} \quad i \neq j \\
	0  & \text{if} \quad i=j,
\end{cases} 
\end{align}
which mathematically represents the interactions in the vortical network.  Since a vortex cannot induce motion upon itself, there are no diagonal entries for the adjacency matrix. Here we have an undirected network with a symmetric adjacency matrix $\mathbf{A}$ to enable the use of the algorithm described later. However, it is possible to define an asymmetric adjacency matrix for a directed network.

Given the vorticity field of a flow, the corresponding vortical interaction network can be formulated as above. This also suggests that the size of the adjacency matrix scales as $\mathbf{A} \in \mathbb{R}^{n \times n}$ where $n$ is the number of grid points or vortical elements in the flow field. For a high-fidelity simulation or experimental measurements of a flow field, $n$ can exceed $\mathcal{O}(10^5-10^6)$, making the computation and analysis of the vortical network computationally expensive. Below, we discuss a few network measures that are particularly useful for understanding the interaction-based characteristics of the flow. The current work focuses on computing these measures for large, dense vortical networks in a computationally tractable manner.}

\subsection{Network measures and computations}
	
%
The adjacency matrix is the basis for a number of useful network measures, such as eigenvector centrality, Katz centrality, and PageRank \cite{Newman10,page1999pagerank}. {The centrality of a network describes the most important or central elements based on measures that quantify the connection strength  of the nodes. One of the most basic centrality measures is the node strength, given by $s_i = \sum_j A_{ij}$, which measures the overall interaction strength of a node. The eigenvector centrality is another fundamental network property, which also considers the strength of the neighboring elements of a node and is used to identify the most influential nodes in the network \cite{kolaczyk2014statistical}. It is computed as the eigenvector of $\mathbf{A}$ corresponding to the largest eigenvalue $\lambda$:
\begin{align}
\mathbf{A}\mathbf{u} = \lambda\mathbf{u}.
\end{align}
Entries with large absolute value in the vector $\mathbf{u}$  correspond to influential nodes.}
The Katz centrality and PageRank are related and also based on computing the dominant eigenvector of the adjacency matrix.

{Network analysis is also useful for graph partitioning and clustering \cite{Newman10}, where the nodes in the network are grouped into clusters based on their close connectivity. Spectral partitioning or clustering is a common approach based on the dominant eigenvector of the adjacency matrix. Nodes can be divided and sub-divided into groups based on the sign of the corresponding elements in the dominant eigenvectors of $\mathbf{A}$. Spectral partitioning can also be used in community detection algorithms to identify groups of closely interacting nodes, called communities \cite{clauset2004finding}. 
These algorithms perform a hierarchical grouping of nodes into communities based on maximizing the overall modular measure of the network.
This hierarchical grouping can be performed using the dominant eigenvectors of $\mathbf{A}$ \cite{leicht2008community}. In vortical flows, these network-based measures are useful for identifying the most influential vortical elements and clustering them into coherent structures or communities.

Due to the importance of the spectral information of the adjacency matrix, considerable effort has gone into developing accurate and efficient algorithms to compute dominant eigenvalues for large adjacency matrices.} 
%
Provided that the largest eigenvalue is real and distinct, the dominant eigenvector can be efficiently computed using the power method~\cite{trefethen1997numerical}.  
A symmetric adjacency matrix with nonnegative entries will have purely real eigenvalues. 
For instance, the power method is commonly used to compute the dominant eigenvector of an adjacency matrix~\cite{bryan200625}, and this method is particularly efficient for large sparse adjacency matrices.
The computational costs of the power method for computing the dominant eigenvector scales as $\mathcal{O}(n^2)$.



\section{Random sampling-based matrix approximation}
{
As stated above, we are interested in computing the eigendecomposition of a real symmetric matrix $\mathbf{A} \in \mathbb{R}^{n\times n}$ that takes the form
\begin{equation}
\mathbf{A} = \mathbf{U} \mathbf{D} \mathbf{U}^\top,
\end{equation}
where the columns of the orthogonal matrix $\mathbf{U} \in \mathbb{R}^{n \times n}$ are eigenvectors and the entries of the diagonal matrix $\mathbf{D} \in \mathbb{R}^{n \times n}$ are the corresponding eigenvalues $\lambda_1 \geq \dots \geq \lambda_{n}$. More concretely, we are interested in computing only the dominant $k$ eigenvectors and eigenvalues (where $k \ll n$) that yield the rank-$k$ approximation
\begin{equation}
\mathbf{A}_k = \mathbf{U}_k  \mathbf{D}_k  \mathbf{U}_k ^\top,
\end{equation}
where $\mathbf{U}_k \in \mathbb{R}^{n \times k}$ contains only the $k$ columns of $\mathbf{U}$ that correspond to the $k$ largest eigenvalues.
%
%
%
%
Here, we use tools from the field of \emph{randomized numerical linear algebra}~\cite{Mahoney2011,Drineas:2016:RRN:2942427.2842602,halko2011finding,Liberty2007pnas,sarlos2006improved,Martinsson201147,erichson_jss,benjamin2017compressed,kannan2017randomized,woodruff2014sketching,woolfe2008fast} to compute such a low-rank approximation efficiently.
Indeed, randomized methods are widely used for computing low-rank approximations~\cite{liberty2013simple,erichson2016randomizedDMD,erichson2017randomizedCP,erichson2018randomizedNMF,erichson2016compressed,shabat2016randomized,rokhlin2009randomized}.

Scalability is achieved by forming a sketched representation of the input (adjacency) matrix which extracts the inherent spectral information.
A sketch $\mathbf{Y}  \in \mathbb{R}^{n \times k}$ can be constructed by post-multiplying the adjacency matrix with a an arbitrary `sketching' matrix $\mathbf{S}  \in \mathbb{R}^{n \times k}$.
The random matrix $\mathbf{S}$ represents either a random projection or a random sampling process.
\begin{itemize}
    \item \emph{Projection-based methods} construct a `sketch' by forming a set of $k$  randomly weighted linear combinations of the columns of the input (adjacency) matrix~\cite{frieze2004fast,woodruff2014sketching,Mahoney2011,Drineas:2016:RRN:2942427.2842602}. This approach can substantially reduce the computational demands when computing a low-rank approximation.  However, projection-based methods generally require at least a single pass over the entire data matrix.
    
    \item \emph{Sampling-based methods} aim to approximate the low-rank structure of the input (adjacency) matrix from a small random subsample of actual columns or rows of the matrix~\cite{ma2015statistical,drineas2006sampling,drineas2012fast}. The sampling process is described in more detail below. In practice, sampling-based methods may bypass the construction of the full adjacency matrix, while sacrificing some accuracy for improved scalability. 
\end{itemize}

In what follows, we generate qualitatively accurate estimates of dominant flow structures for cases where the adjacency matrix is exceedingly large. 
In this problem context, we require a qualitatively correct view of the dominant graph structures, which can be used for the downstream tasks of sensor and actuator placement, optimization, and control.  
Therefore, we investigate the use of column sampling to compute a sketched singular value decomposition and the Nystr\"om Method to approximate the leading $k$ eigenvalues and eigenvectors. 
Both techniques have been successfully used to approximate large-scale matrices in other applications~\cite{talwalkar2013large}. 
It is important to note that sampling-based methods provide tunable error bounds that are based on the number of sampled columns/rows, making them viable as an alternative to well-established deterministic algorithms. 
}

\subsection{Sketched singular value decomposition}
\begin{figure}[t]
\hspace{0.23in}
	\centering
	\begin{subfigure}[t]{0.47\textwidth}
		\centering
		\DeclareGraphicsExtensions{.pdf}
		\hspace{-.6in}\includegraphics[width=.9\textwidth]{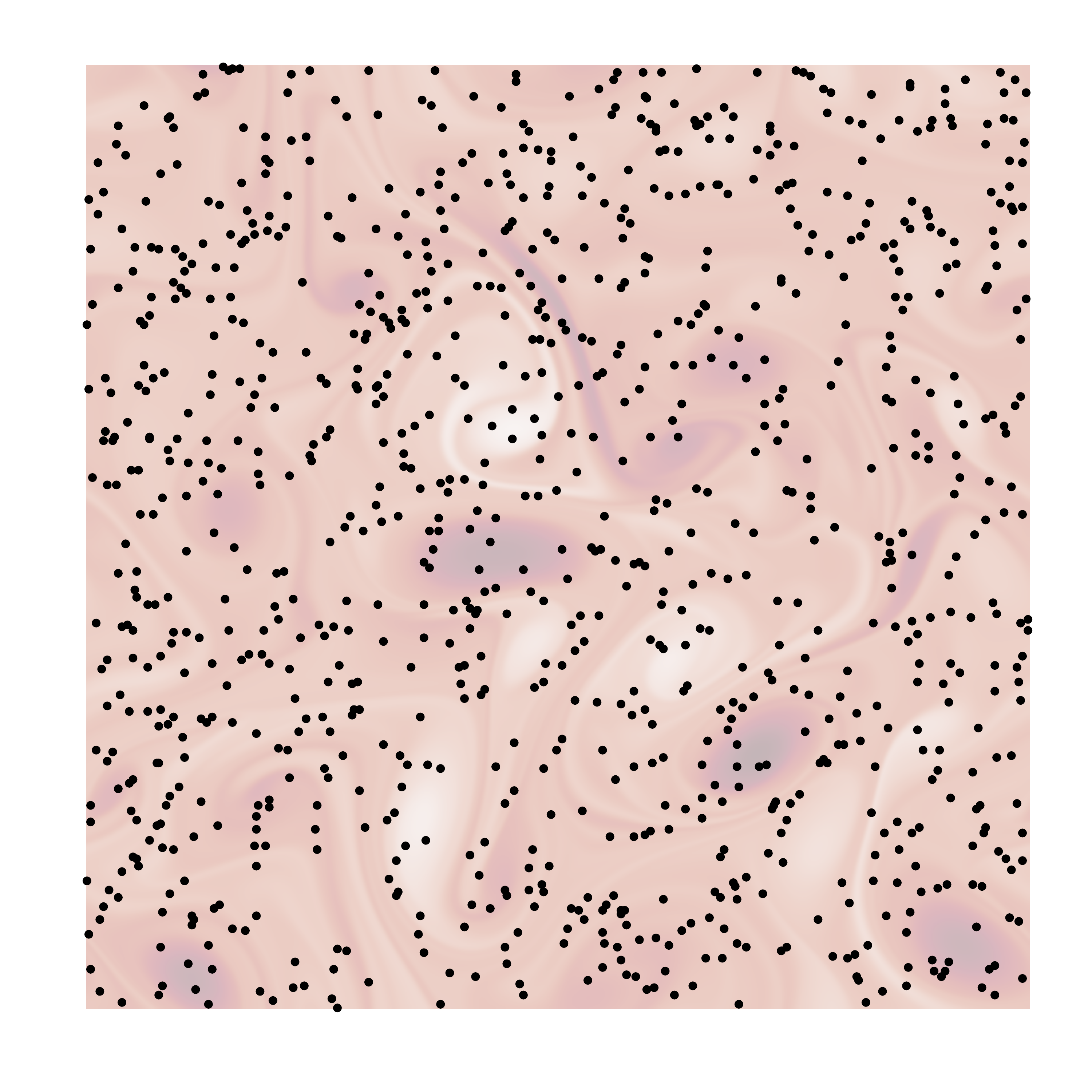}
		\vspace{-.1in}
		\caption{Uniform random sampling. }
	\end{subfigure}
	~
	\begin{subfigure}[t]{0.47\textwidth}
		\centering
		\DeclareGraphicsExtensions{.pdf}
		\hspace{-.6in}\includegraphics[width=.9\textwidth]{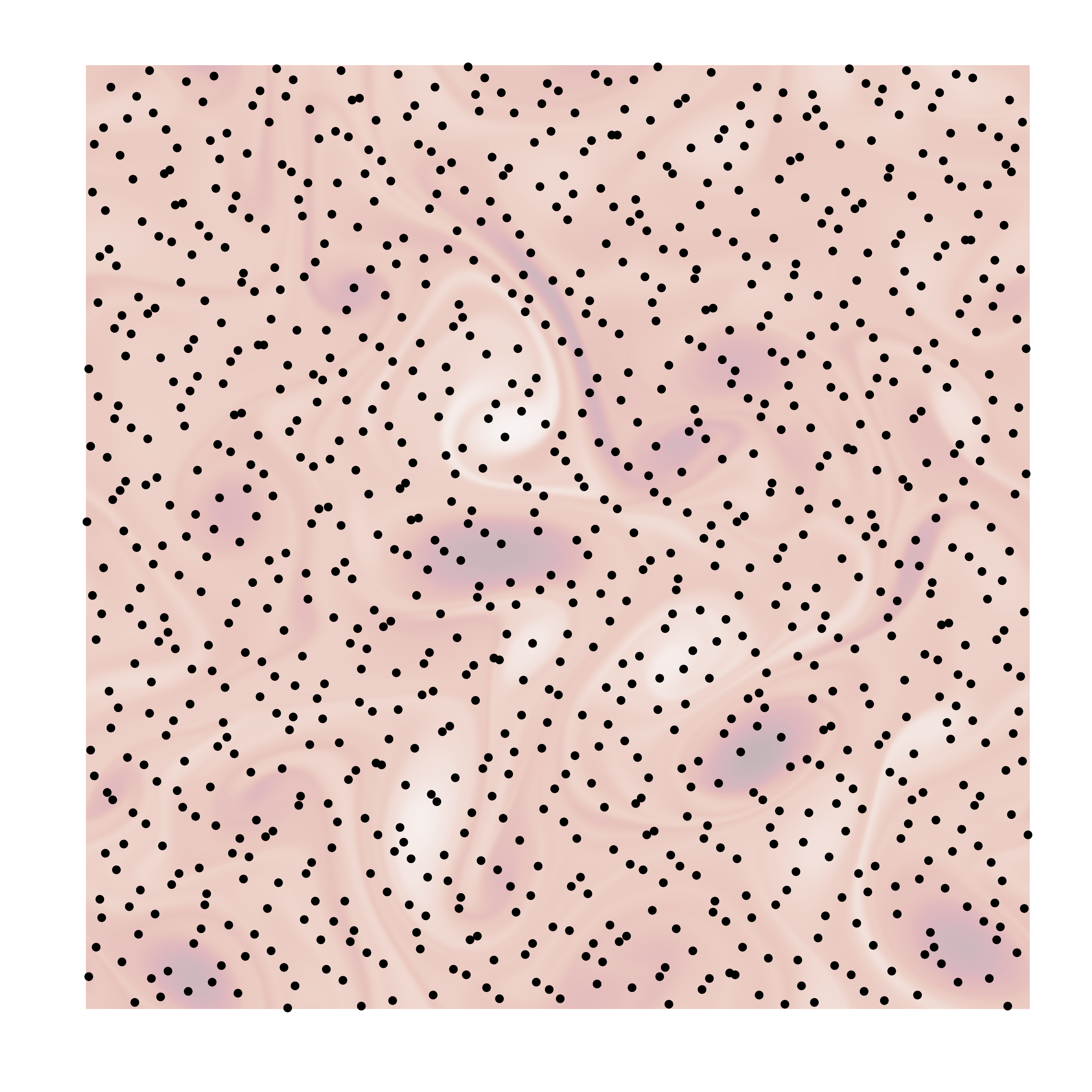}
		\vspace{-.1in}
		\caption{Halton random sampling. }
	\end{subfigure}		
	\vspace{0in}	
	\caption{Examples of two different random sampling approaches based on uniform sampling and Halton random sampling. Halton sampling provides a better coverage of the domain than uniform sampling does.}
	\label{fig:random_numbers}
\end{figure}

{
Column sampling for low-rank matrix approximation dates back to the pioneering work of Frieze et al.~\cite{frieze2004fast}. 
Instead of computing the full eigendecomposition of $\mathbf{A}$, it is possible to form a rank-$k$ approximation by first sampling $l\geq k$ columns from the input matrix 
\begin{equation}
\mathbf{C} = \mathbf{A}\mathbf{S}
\end{equation}
where the matrix $\mathbf{C} \in \mathbb{R}^{n \times l}$ consists of a subset of $l$ columns of $\mathbf{A}$ and $\mathbf{S}$ describes the corresponding random sampling process. In practice, it is only necessary to form a list $J$ comprised of $l$ column indices
\begin{equation}
\mathbf{C} := \mathbf{A}(:,J).
\end{equation}
%
There are several sampling strategies to choose \emph{good} columns, and the approximation error depends on the both the sampling process and the number of columns sampled.  

In our problem setup, each column of $\mathbf{A}$ corresponds to the network interaction of a single grid element in the original fluid flow field with every other grid element. 
Here, we compare uniform random sampling against (quasi-random) Halton sampling~\cite{niederreiter1992random,Halton}. 
Halton sampling provides better coverage of the domain, as illustrated in Figure~\ref{fig:random_numbers}, and our empirical results show that this translates into faster convergence of the approximation error with respect to the number of columns sampled.

If the columns of $\mathbf{C}$ are carefully chosen, then $\mathbf{C}$ provides a basis for the column space of the adjacency matrix $\mathbf{A}$. Note that here we assume that $\mathbf{A}$ is a symmetric matrix.  Thus, the dominant $k$ left singular vectors of $\mathbf{C}=\mathbf{U_C} \boldsymbol{\Sigma}_\mathbf{C} \mathbf{V^*_C}$, provide an approximation for the dominant $k$ eigenvectors of $\mathbf{A}$, i.e., we have $\tilde{\mathbf{U}}_k \approx \mathbf{U_C}[:,1:k]$. 
The corresponding eigenvalues are approximated by scaling the singular values of $\mathbf{C}$:
\begin{align}
\tilde{\mathbf{D}}_k= \sqrt{\frac{n}{l}}\boldsymbol{\Sigma}_{\mathbf{C}}[1:k,1:k].
\end{align}
The cost of constructing the sketch matrix using uniform or Halton sampling is $\mathcal{O}(nl)$. Then, the cost of computing the SVD on $\mathbf{C}$ is $\mathcal{O}(nl^2)$.
}
For parallel computations, it is possible to compute $ \mathbf{U_C}$ and $\boldsymbol{\Sigma}_\mathbf{C}$ by taking the SVD of $\mathbf{C}^T\mathbf{C}$, which requires $\mathcal{O}(nl^2)$ to be generated and $\mathcal{O}(l^3)$ for the corresponding SVD.

\subsection{The Nystr\"om method}
The Nystr\"om method provides an efficient approach for low-rank matrix approximations in large-scale learning applications. 
It was initially introduced as a quadrature method for numerical integration to approximate eigenfunction solutions. 
Recent work has streamlined these computations and provided theoretical foundations for the use of various sampling schemes~\cite{williams2001using,drineas2005nystrom,zhang2008improved,gittens2016revisiting}.

In addition to column sampling, the Nystr\"om method further constructs a square matrix $\mathbf{W}\in \mathbb{R}^{l \times l}$ by selecting the $J$ rows and columns of $\mathbf{A}$:
\begin{equation}
\mathbf{W} := \mathbf{C}(J,:)=\mathbf{A}(J,J).
\end{equation}
When $\mathbf{A}$ is symmetric and positive-semidefinite (SPSD), $\mathbf{W}$ is also SPSD. 
The small matrix $\mathbf{W}$ can be used to efficiently compute the dominant eigenvectors and eigenvalues of $\mathbf{A}$. 
Following~\cite{kumar2012sampling}, we first compute the eigendecomposition:
\begin{equation}
\mathbf{W} = \mathbf{{U}}_\mathbf{W}  \mathbf{{D}}_\mathbf{W}   \mathbf{{U}}_\mathbf{W}  ^\top.
\end{equation}
Then, we reconstruct the dominant eigenvalues as
\begin{equation}
\tilde{\mathbf{D}}_k = \frac{n}{l} \mathbf{{D}}_\mathbf{W} ,
\end{equation}
and the corresponding eigenvectors
\begin{equation}
\tilde{\mathbf{U}}_k = \sqrt{\frac{l}{n}} \mathbf{C} \mathbf{{U}}_\mathbf{W}  \mathbf{{D}}^\dagger_\mathbf{W} .
\end{equation}
%
The Nystr\"om method has a computational complexity of $\mathcal{O}(l^3)$.


\section{Results}
We now demonstrate the use of sampling-based randomized decomposition techniques to characterize the vortical interaction networks for two example flows: the laminar wake flow past a NACA 0012 airfoil with a Gurney flap attached to the trailing edge~\cite{gopalakrishnan2017airfoil}, and a two-dimensional homogeneous decaying isotropic turbulence~\cite{Taira2016JFM}.  

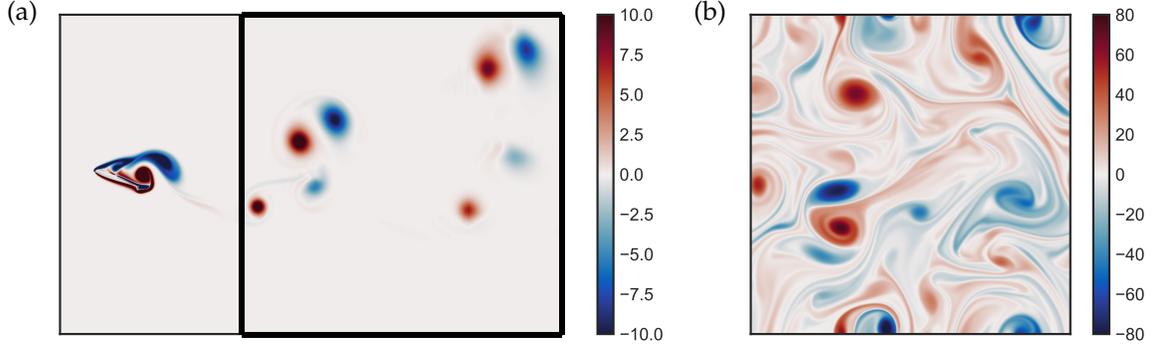
\begin{figure}[t]
	\centering
	\DeclareGraphicsExtensions{.pdf}
	\begin{overpic}[width = 0.6\textwidth]{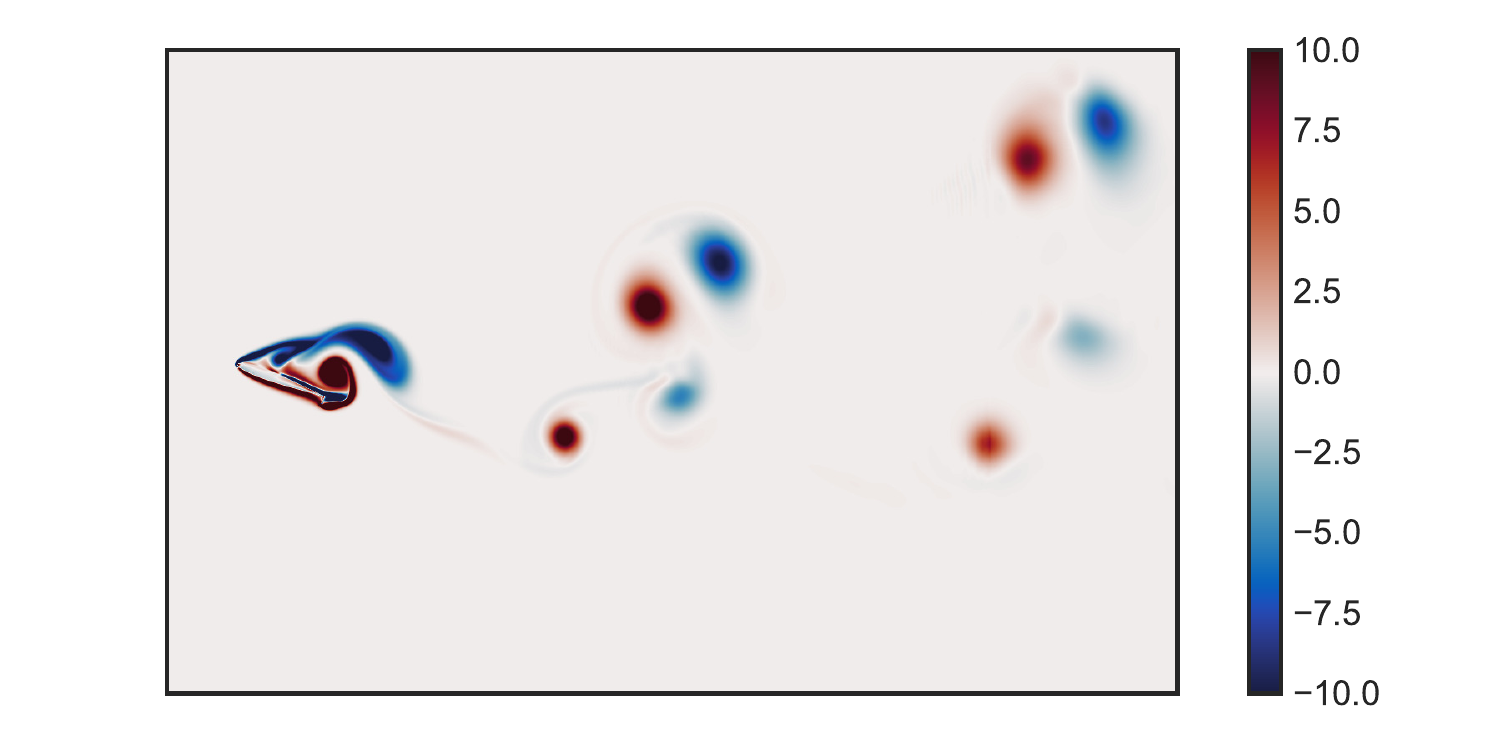}
	\small 
	\put(4,46){(a)}
	\linethickness{2pt}
 	\put(35.5,3.7){\line(1,0){43}}
 	\put(78.5,3.7){\line(0,1){43}}
    \put(78.5,46.7){\line(-1,0){43}}
    \put(35.5,46.7){\line(0,-1){43}}
	\end{overpic}
	\hspace{-0.7cm}
	\begin{overpic}[width = 0.4\textwidth]{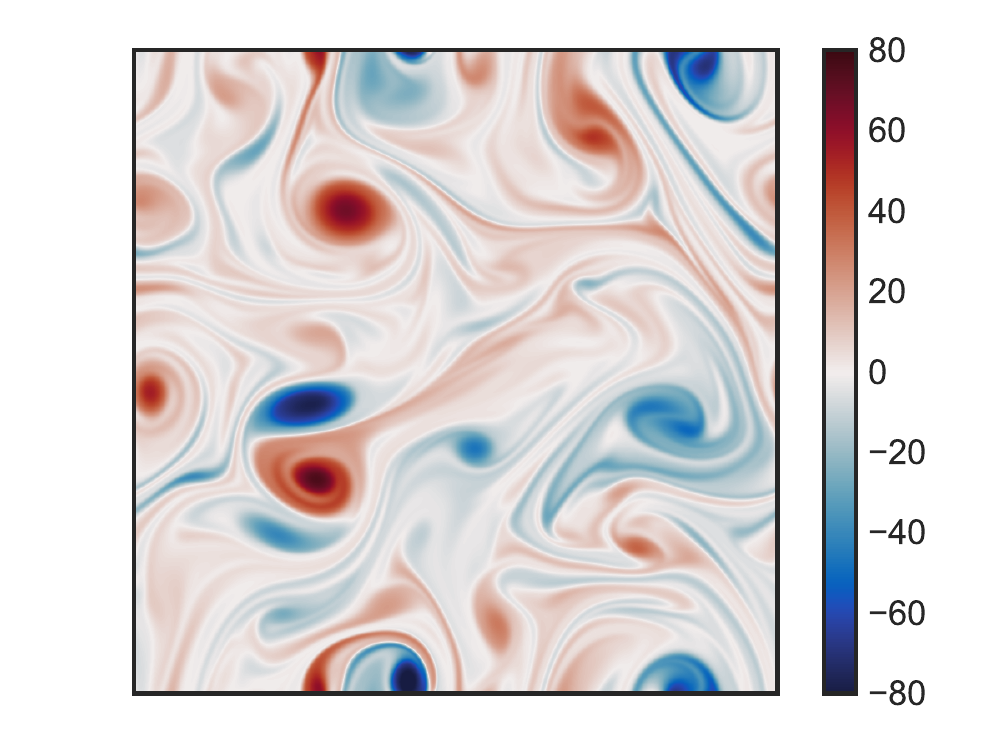}
	\small\put(2,69){(b)}
	\end{overpic}
	\caption{Example flow fields.  (a) Vorticity field of two-dimensional DNS of the flow over a NACA 0012 airfoil with Gurney flap at an angle of attack of $20^\circ$ and flap height of $0.1$ chord length at $Re=1000$; the full illustration can be found in~\cite{gopalakrishnan2017airfoil}. A subdomain of the original vorticity field is used in this example. (b) Vorticity field of two-dimensional decaying homogeneous isotropic turbulence, initialized by an integral length-scale based Reynolds number of $Re(t_0) = 814$.}
	\label{fig:airfoil}
\end{figure}

\subsection{Numerical simulations}\label{eample_2}
The first example is the flow past a NACA 0012 airfoil with a chord-based Reynolds number of $Re = 1,000$, at an angle of attack of $20^{\circ}$, and with a Gurney flap of $10\%$ chord length ($c$).
This setup generates an unsteady 2P wake \cite{williamson1988vortex,gopalakrishnan2017airfoil}, with periodic shedding of two pairs of positive and negative vortices. 
This flow is solved using direct numerical simulations (DNS) via the immersed boundary projection method \cite{taira:07ibfs,taira:fastIBPM}, following the computational setup of Gopalakrishnan Meena \textit{et al.}~\cite{gopalakrishnan2017airfoil}.  
%
The computational domain consists of five nested levels of multidomains with the finest level of size $(x/c,y/c) \in [-1,1] \times [-1,1]$ and the largest being $(x/c,y/c) \in [-16,16] \times [-16,16]$ in size. All the domains have a grid resolution of $m_x \times m_y =  360 \times 360$.
The time step is limited to a maximum Courant-Friedrichs-Lewy (CFL) number of $0.3$.
Figure~\ref{fig:airfoil} (a) shows the vorticity field in a sub-domian with a grid resolution of $250\times 150$ and size $5.53\times 3.31$, nondimensionalized by the chord length of the airfoil.


The second example is two-dimensional decaying homogeneous isotropic turbulence, which is a canonical high-dimensional, mutiscale turbulent flow that exhibits complex nonlinear interaction of vortices over a wide range of length scales \cite{boffetta2012two}. 
The vorticity field is obtained from a two-dimensional incompressible DNS, solving the vorticity transport equation without forcing \cite{Taira2016JFM}. 
The simulation, based on the Fourier spectral method and the fourth-order Runge-Kutta time integration scheme, is performed on a square biperiodic computational domain 
with grid resolution of $m_x \times m_y = 1024 \times 1024$. 
The vorticity field is initialized with $100$ vortices with random strengths, core sizes, and centers such that the kinetic energy spectra satisfies $E(k) \propto k \exp(-k^2/k_0^2)$ with $k_0 = 26.5$ \cite{kida1985numerical}. 
For this study, the flow field is initialized by an integral length-scale based Reynolds number of $Re(t_0) = 814$, shown in Figure~\ref{fig:airfoil} (b). 

{The vorticity snapshots of both flows are used to construct the adjacency matrices of the corresponding vortical networks using Eq.~\ref{e:adj}. The size of the adjacency matrices are of the order $\mathcal{O}(10^8)$ and $\mathcal{O}(10^{12})$ for the airfoil and turbulent flow, respectively, emphasizing the need for dimensionality reduction and computationally tractable analysis tools.

\subsection{Computational performance}

A naive application of network theory in fluid dynamics can be challenging, because the adjacency matrix is often too large to be stored or analyzed using classical techniques from linear algebra. 
There are two primary challenges associated with large vortical networks:
\begin{itemize}
	\item The size of the adjacency matrix, which describes the network edges, scales as the square of the number of fluid grid points;	
	\item Unlike social graphs, the vortical adjacency matrix is dense.
\end{itemize}
\begin{figure}[!b]
	\centering	
	\begin{subfigure}[t]{1\textwidth}
		\centering
		\DeclareGraphicsExtensions{.pdf}
		\includegraphics[width=0.47\textwidth]{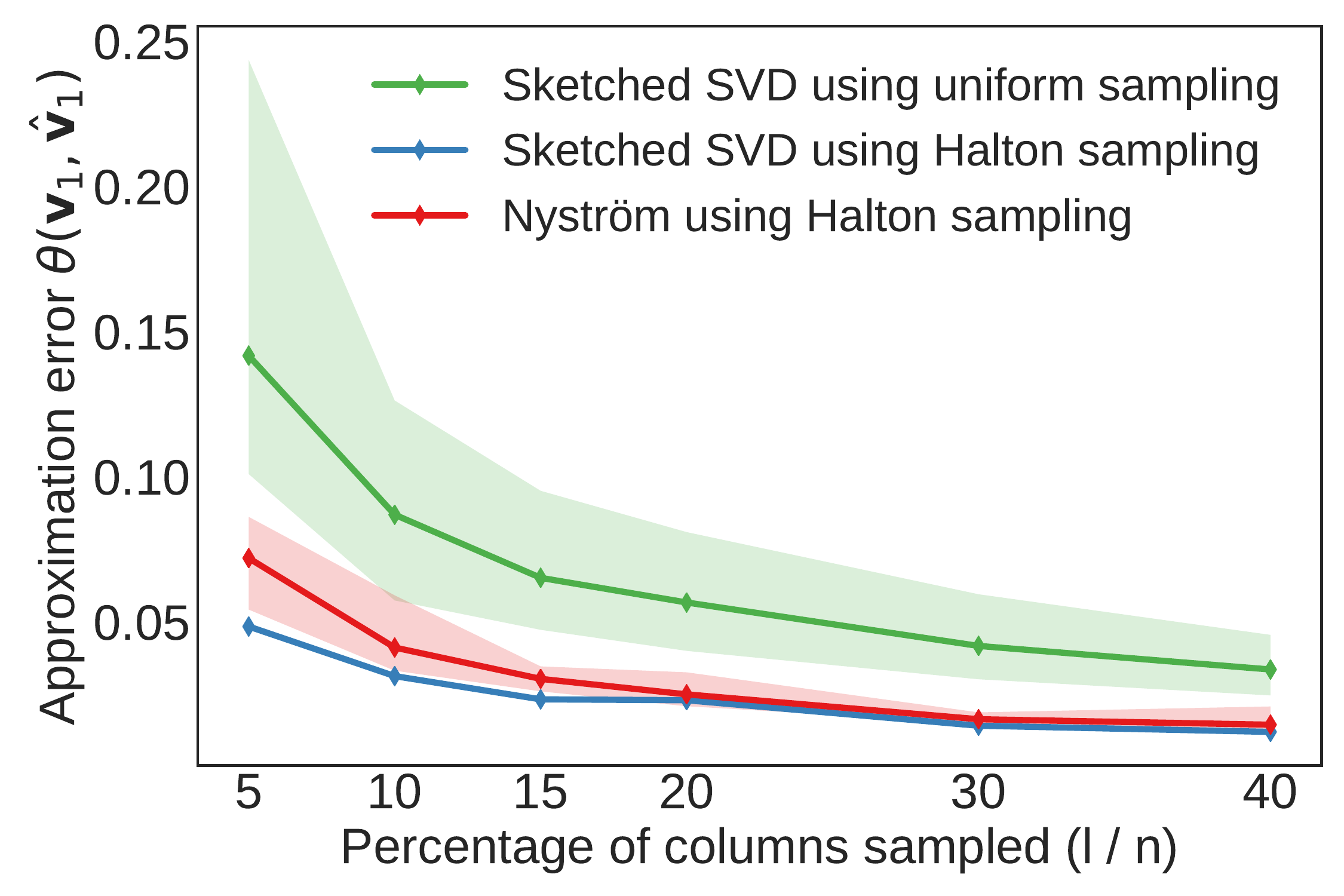}
		\includegraphics[width=0.47\textwidth]{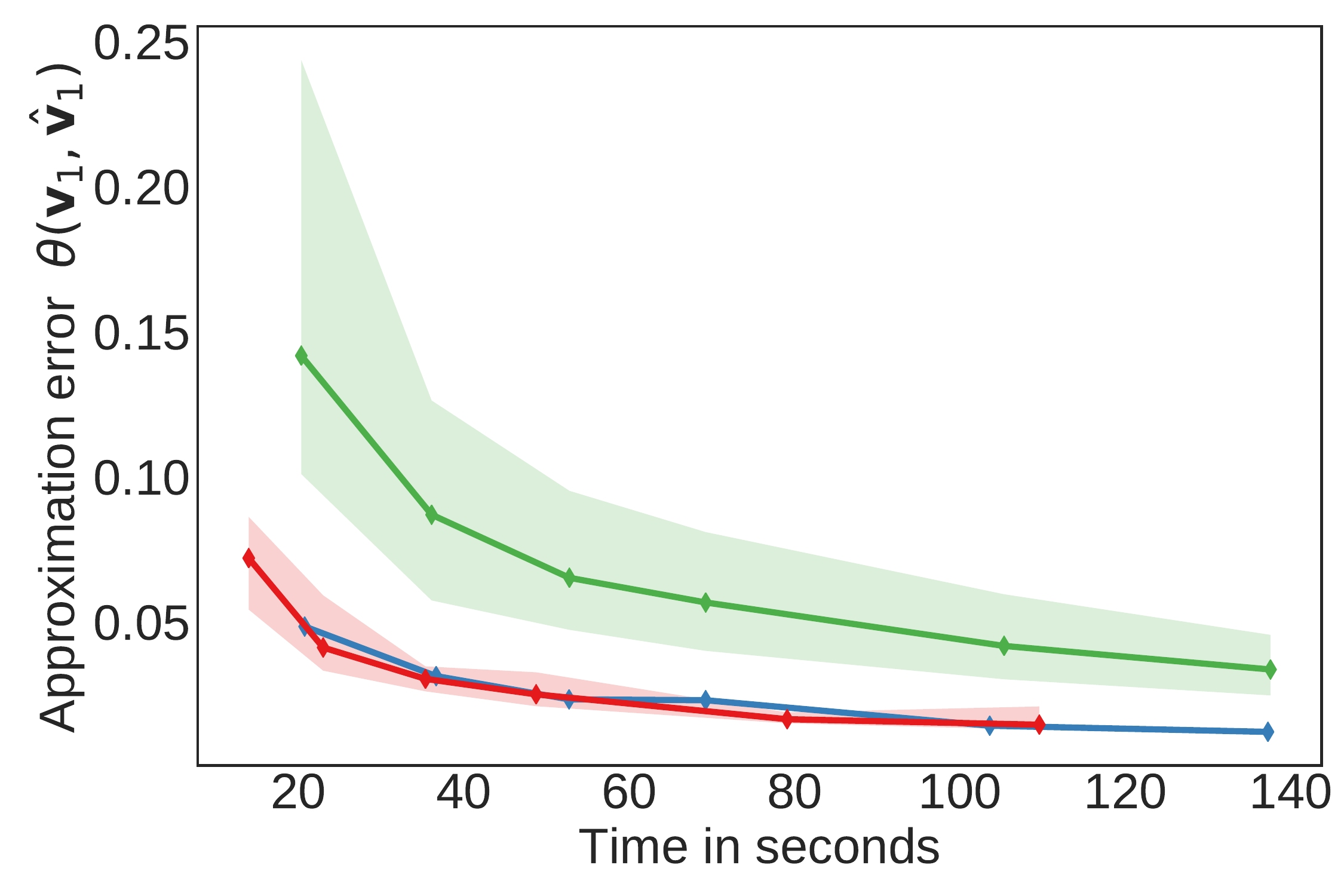}
		\vspace{-.05in}\caption{Results for an $n\times n$ adjacency matrix with $n=316^2$. }
		\label{fig:simulation_A_error}
	\end{subfigure}		
	\begin{subfigure}[t]{1\textwidth}
		\centering
		\DeclareGraphicsExtensions{.pdf}
		\includegraphics[width=0.47\textwidth]{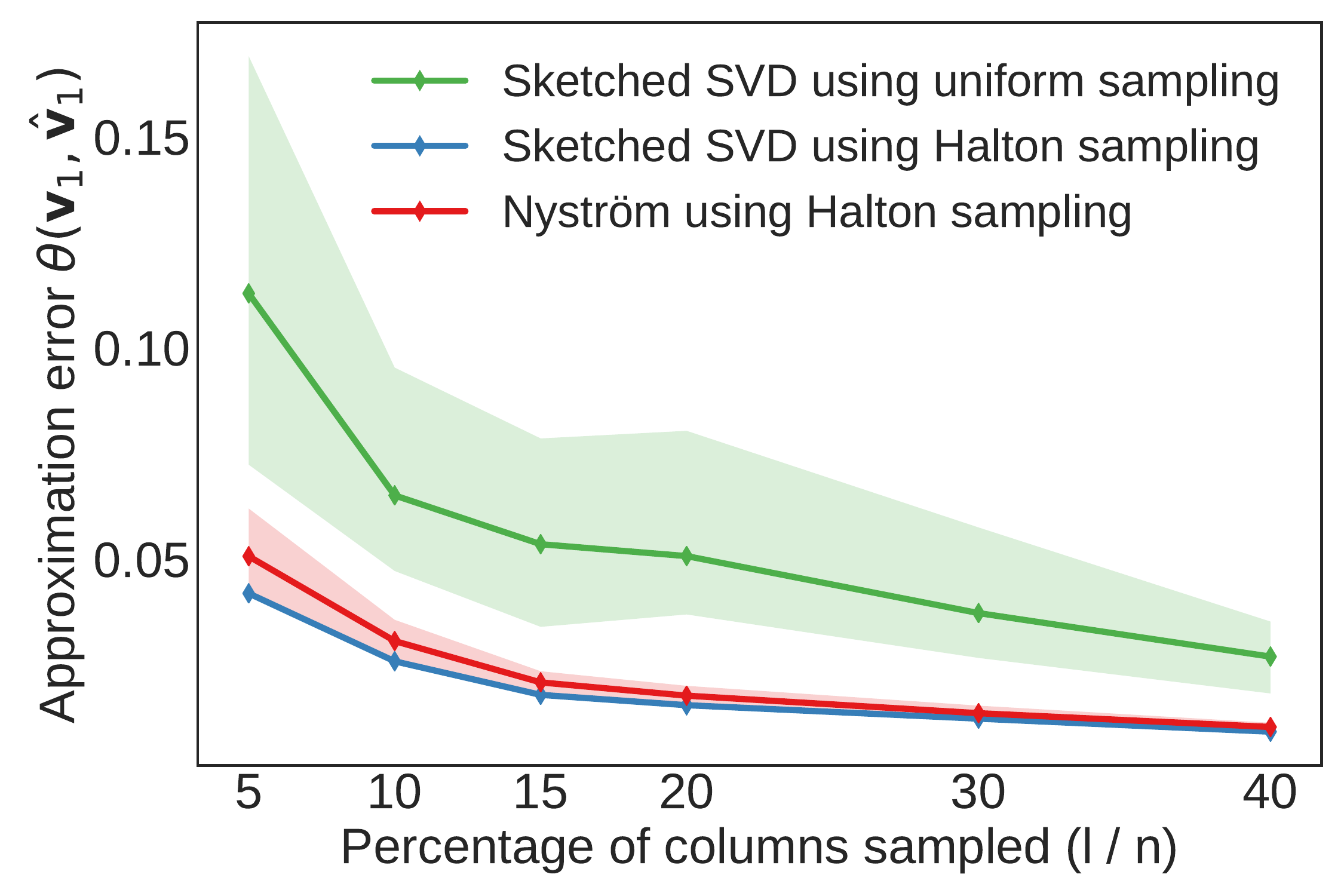}
		\includegraphics[width=0.47\textwidth]{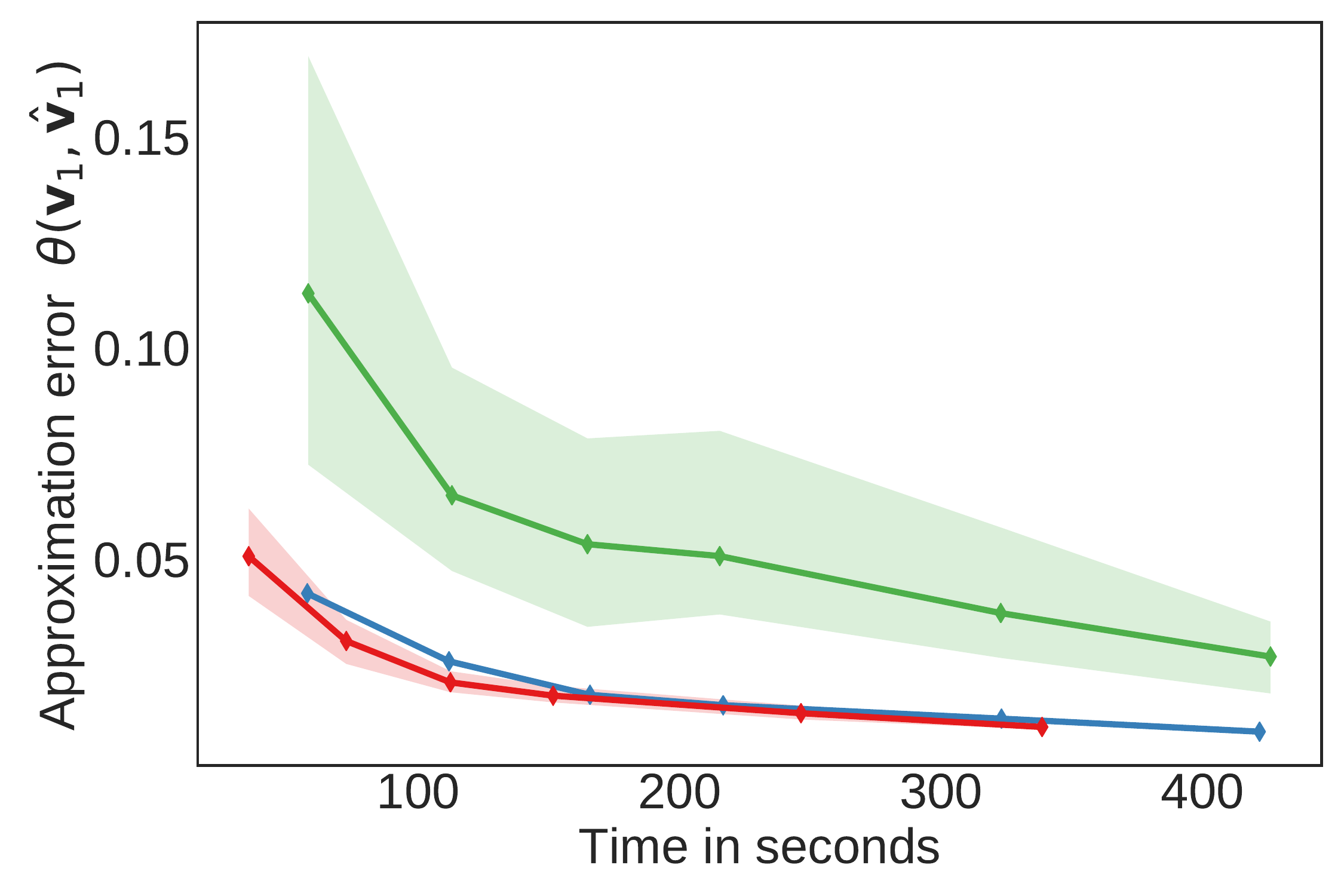}
		\vspace{-.05in}\caption{Results for an $n\times n$ adjacency matrix with $n=421^2$. }
		\label{fig:simulation_A_error_large}
	\end{subfigure}
	\vspace{0.05in}
	\caption{Computational performance for the flow over an airfoil using two different spatial resolutions. The results are averaged over $20$ runs using different random seeds.}
	\label{fig:simulation_A}
\end{figure}

Here, we use randomized methods to accelerate computations based on the vortical adjacency matrix.
We sample a small number of columns from the adjacency matrix to compute a sketched SVD or Nystr\"om approximation. 
Thus, the approximation quality can be controlled via the number of samples used to compute the low-rank approximation. 
In the context of networks arising in fluid mechanics, it turns out that sampling a small number of columns is sufficient for computing an approximation that is qualitatively useful. 
%
%
This leads to tremendous memory and computational savings, since the costs are proportional to the number of columns that we sample. 
We run all of our experiments using Amazon AWS, working on a memory optimized instance `x1.16xlarge' that is powered by two Intel Xeon E7 8880 v3 (Haswell) processors with $64$ virtual cores and $976$ GB fast memory. Our algorithms are implemented in Python powered by MKL accelerated linear algebra routines, and C extensions for constructing the elements of the adjacency matrix.

We show in Figure~\ref{fig:simulation_A} the approximation error as a function of the percentage of columns sampled and the corresponding computational time required. 
Here, we consider two cropped spatial grids of different dimensions from which we construct the adjacency matrix. 
The error is reported as the acute angle between the true leading eigenvector ${\mathbf{v}}$ and the approximate eigenvector $\hat{\mathbf{v}}$, which provides a useful summary of how closely the two high-dimensional eigenvectors align.
Since the algorithms are fundamentally based on random sampling, we average the results over $20$ initializations and report the distribution of errors. 

Halton sampling results in considerably lower error than uniform sampling at the same percentage of columns sampled, which is consistent with the observation that it is sampling the flow domain more efficiently.
The results based on Halton sampling also have considerably lower variance, indicating improved numerical robustness.
Figures~\ref{fig:simulation_A_error} and~\ref{fig:simulation_A_error_large} illustrate how substantially the computational time increases when the spatial resolution of the flow field is increased only slightly.
In the $421\times 421$ case, only about $5\%$ of the columns need to be sampled for $5\%$ error compared between the approximate and true eigenvectors.
In other words, we use only about $5\%$ of the information to compute the low-rank approximation, and this reduces the memory requirements for the example in Figure~\ref{fig:simulation_A_error_large} from about $215$GB to roughly $10.5$GB. 
At the same time, the computational time is reduced since it is not necessary to construct the full adjacency matrix. 
The computational time is reduced further by using the Nystr\"om method, while only sacrificing a small amount of accuracy. 
For comparison, Table~\ref{tab:summary} summarizes the computational demands and performance for the deterministic power iteration method.
Figure~\ref{fig:simulation_B} shows similar results for the turbulent flow. Again, we see that Halton sampling is favorable and that only a small fraction of columns is required to yield a good approximation of the leading eigenvector.

}

\begin{figure}[!b]
	\centering	
	\begin{subfigure}[t]{1\textwidth}
		\centering
		\DeclareGraphicsExtensions{.pdf}
		\includegraphics[width=0.47\textwidth]{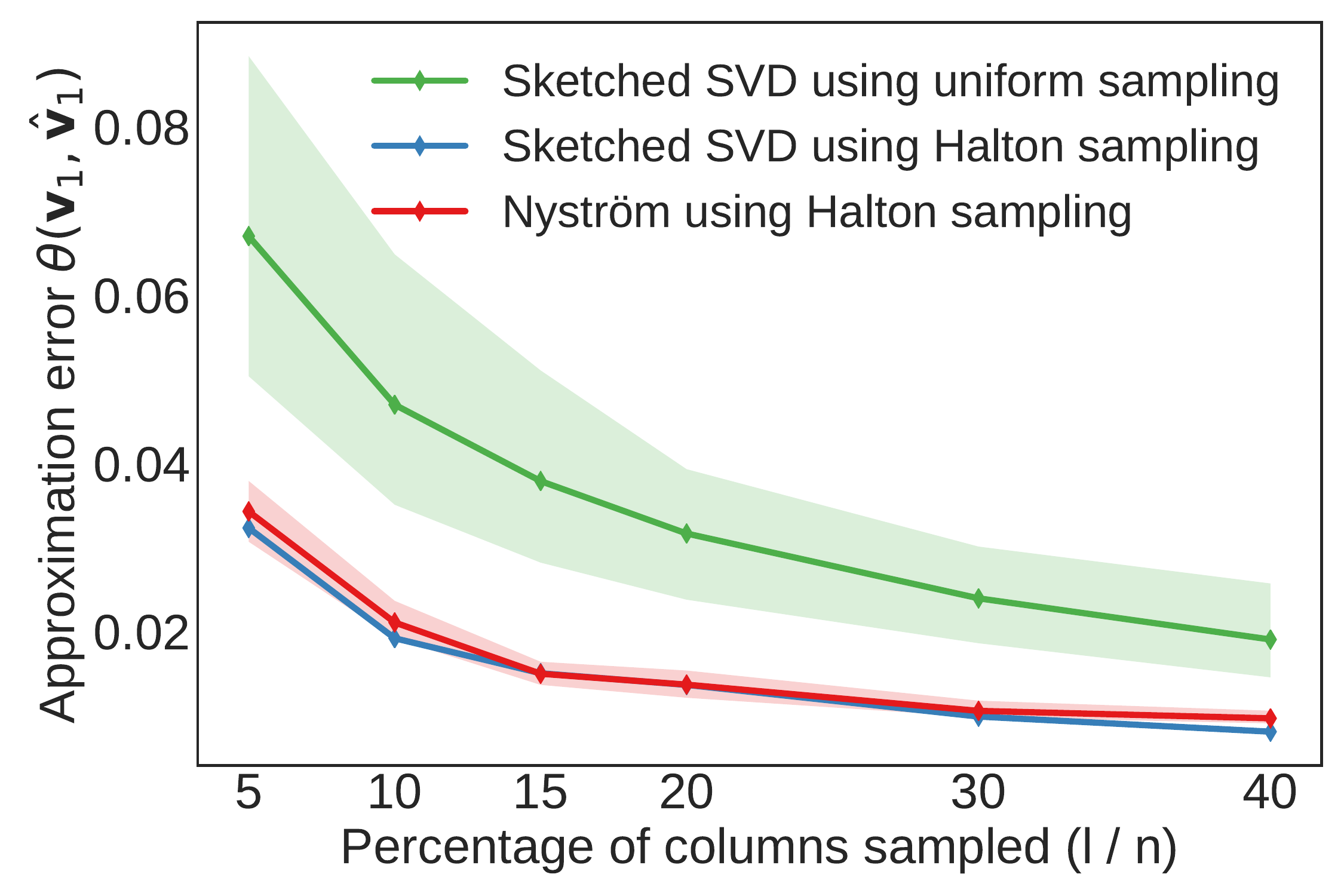}
		\includegraphics[width=0.47\textwidth]{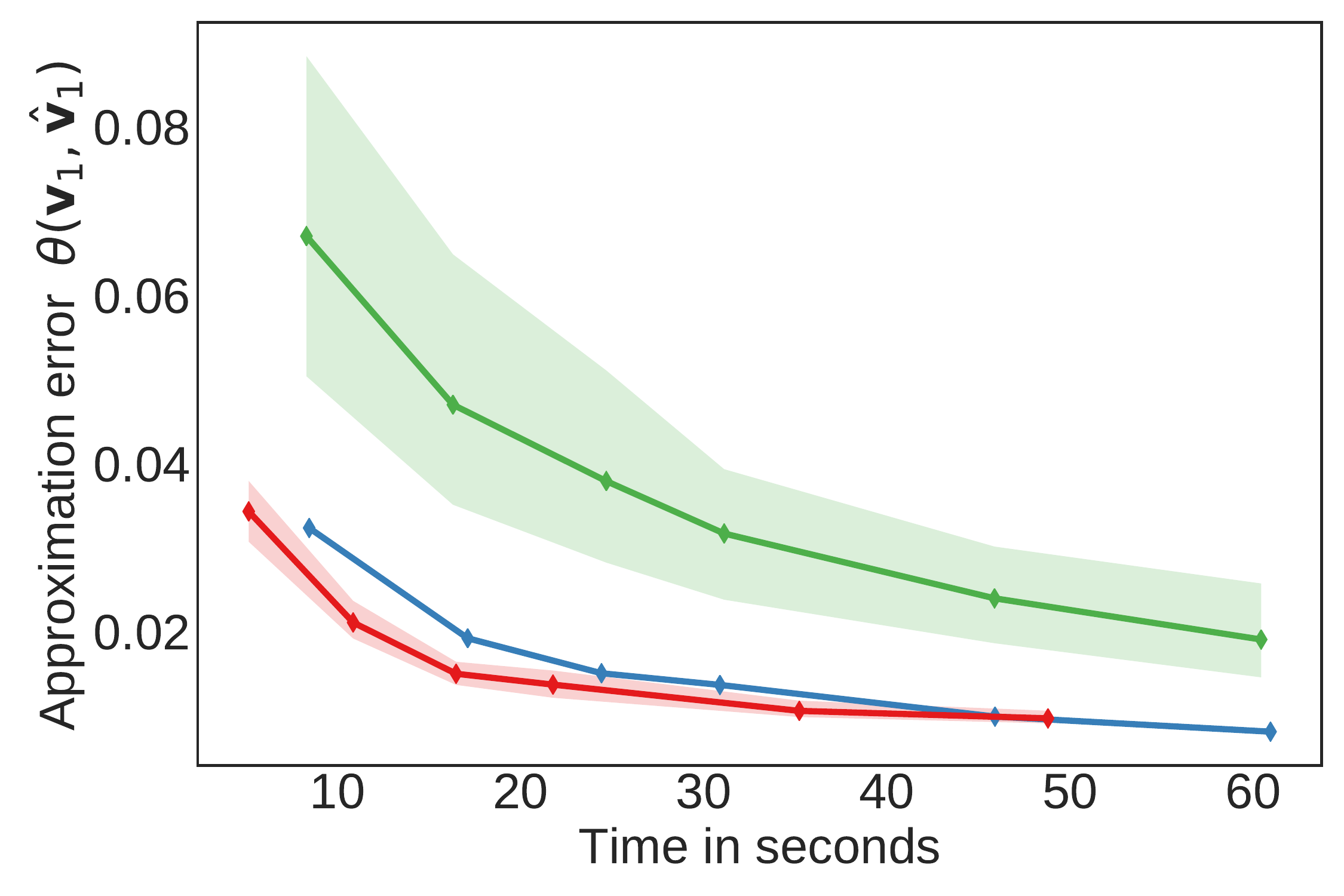}
		\caption{For an $n\times n$ adjacency matrix with $n=256^2$.}
		\label{fig:simulation_B_error}
	\end{subfigure}		
	\begin{subfigure}[t]{1\textwidth}
		\centering
		\DeclareGraphicsExtensions{.pdf}
		\includegraphics[width=0.47\textwidth]{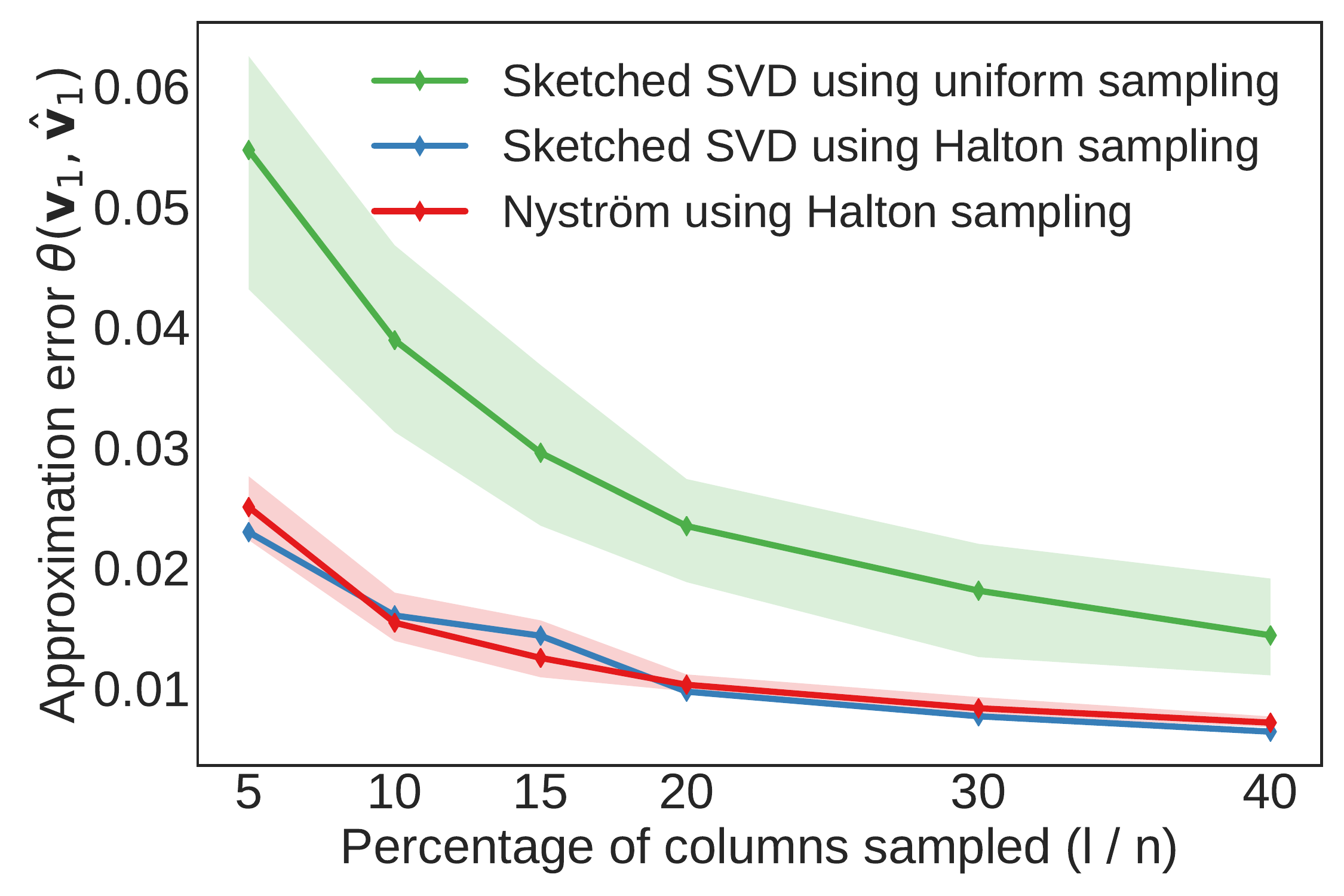}
		\includegraphics[width=0.47\textwidth]{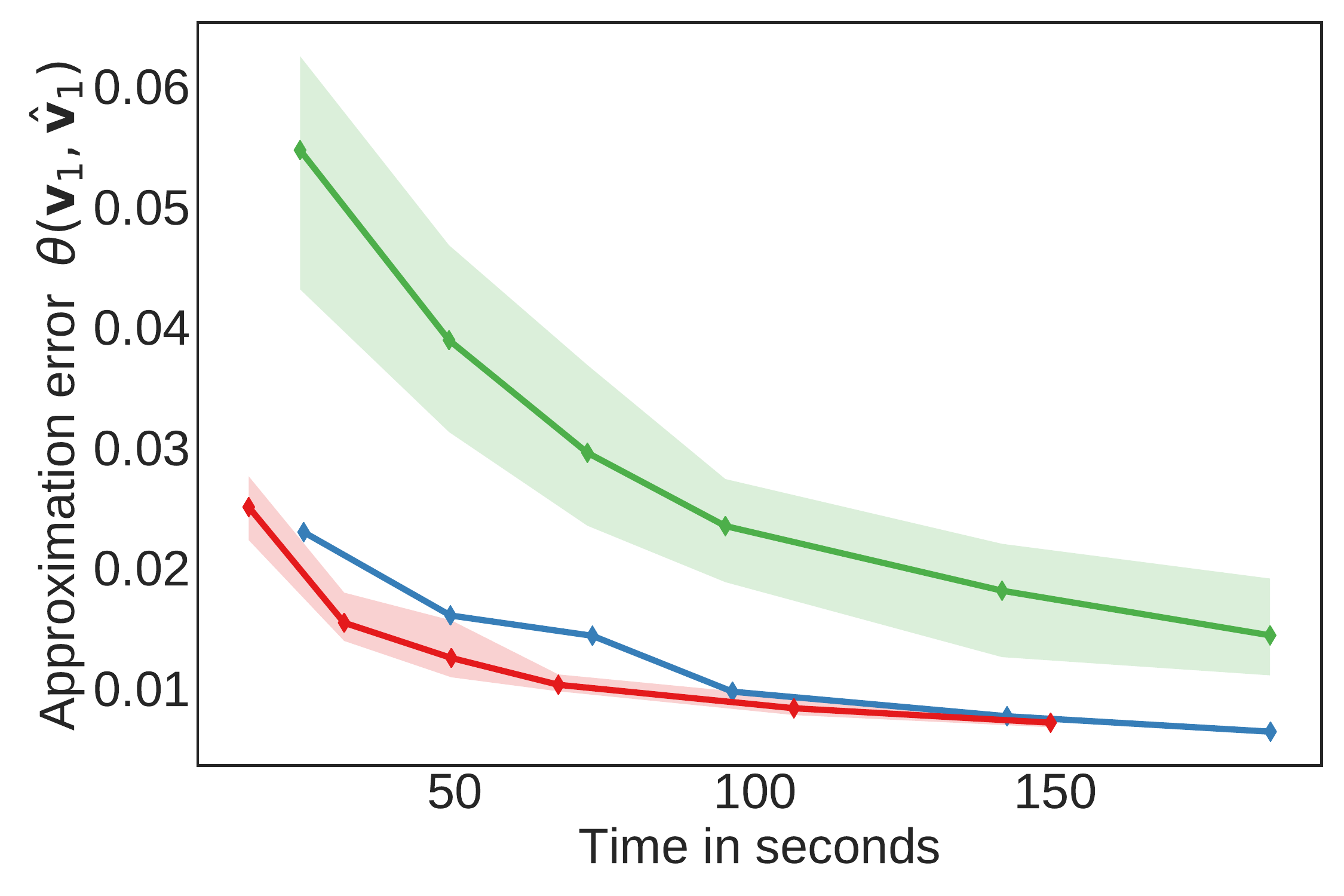}
		\caption{For an $n\times n$ adjacency matrix with $n=342^2$.}
		\label{fig:simulation_B_error_large}
	\end{subfigure}
	\vspace{.05in}	
	\caption{Computational performance for the isotropic turbulent flow using two different spatial resolutions. The results are averaged over $20$ runs using different random seeds.}
	\label{fig:simulation_B}
\end{figure}

\begin{table}[!t]
	\centering\scalebox{0.8}{	
		\begin{tabular}{l ccccc} \toprule
			              & \makecell{Dimensions of \\ adjacency matrix}  & \makecell{Time to construct\\ full adjacency matrix}  & \makecell{Storage \\ required}  & \makecell{Time to compute \\deterministic eigenvector}\\
			\midrule
			 airfoil   & $316^2 \times 316^2$   &  $197.76$ (sec)   & 74.29 (GB) &  $60.84$ (sec)  \\
			 airfoil   & $421^2 \times 421^2$   &  $715.78$ (sec)   & 214.67 (GB) &  $197.87$ (sec) \\
			\midrule
			 turbulent   & $256^2 \times 256^2$   &  $158.44$ (sec)   & 32.00 (GB) &  $27.51$ (sec) \\
			 turbulent   & $342^2 \times 342^2$   &  $272.09$ (sec)   & 101.92 (GB) &  $85.01$ (sec) \\ \bottomrule 
	\end{tabular}}
	\vspace{+.10in}
	\caption{Summary of computational results for constructing the full adjacency matrix and computing the dominant eigenvector using the deterministic power method. Here the computational bottleneck is the memory required to construct the adjacency matrix. Note that the power method does not generally require that the full adjacency matrix is constructed explicitly; however, the power method is not efficient because the adjacency matrices we consider are dense.}
	\label{tab:summary}
\end{table}

\subsection{Approximate eigenvectors and spectral clustering}

\begin{figure}[!b]
    \vspace{-.15in}
    \hspace{.1in}
	\centering	
	\begin{subfigure}[t]{0.315\textwidth}
		\centering
		\DeclareGraphicsExtensions{.pdf}
		 \hspace{-.1in}\tikz[baseline=(a.north)]\node[yscale=-1,inner sep=0,outer sep=0](a)
		{\begin{overpic}[width=1\textwidth]{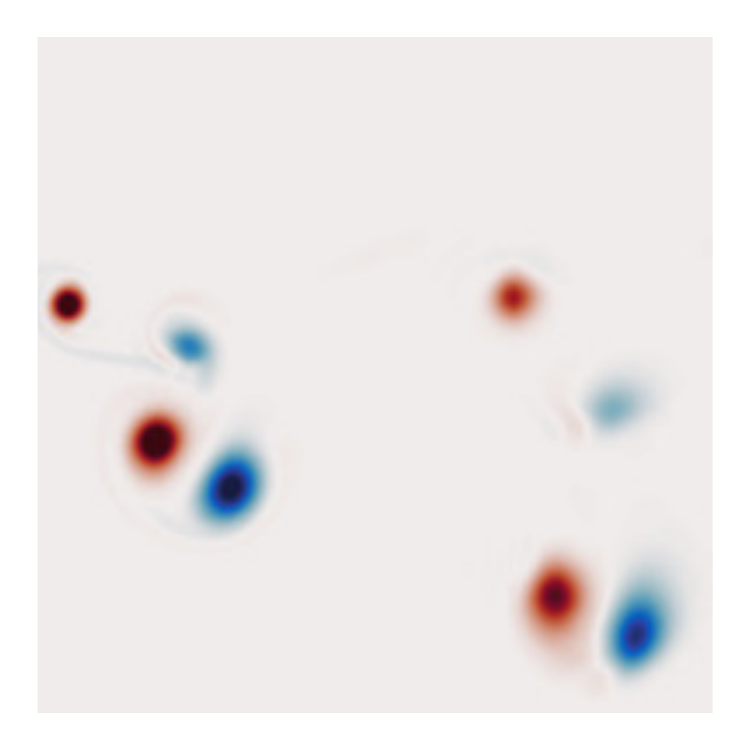}
		\put(-5,49){\reflectbox{\rotatebox[origin=c]{270}{Airfoil wake}}}
		\end{overpic}};
	\end{subfigure}
	~
	\begin{subfigure}[t]{0.315\textwidth}
		\centering
		\DeclareGraphicsExtensions{.pdf}
		 \hspace{-.1in}\tikz[baseline=(a.north)]\node[yscale=-1,inner sep=0,outer sep=0](a)
		{\begin{overpic}[width=1\textwidth]{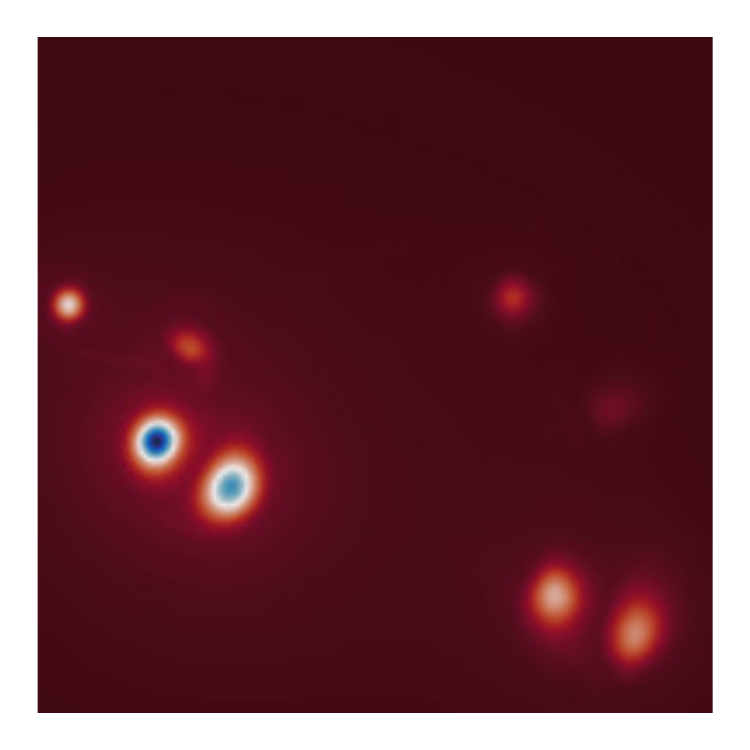}
		\end{overpic}};
	\end{subfigure}
	~ 
	\begin{subfigure}[t]{0.315\textwidth}
		\centering
		\DeclareGraphicsExtensions{.pdf}
		\tikz[baseline=(a.north)]\node[yscale=-1,inner sep=0,outer sep=0](a)
		{\includegraphics[width=1\textwidth]{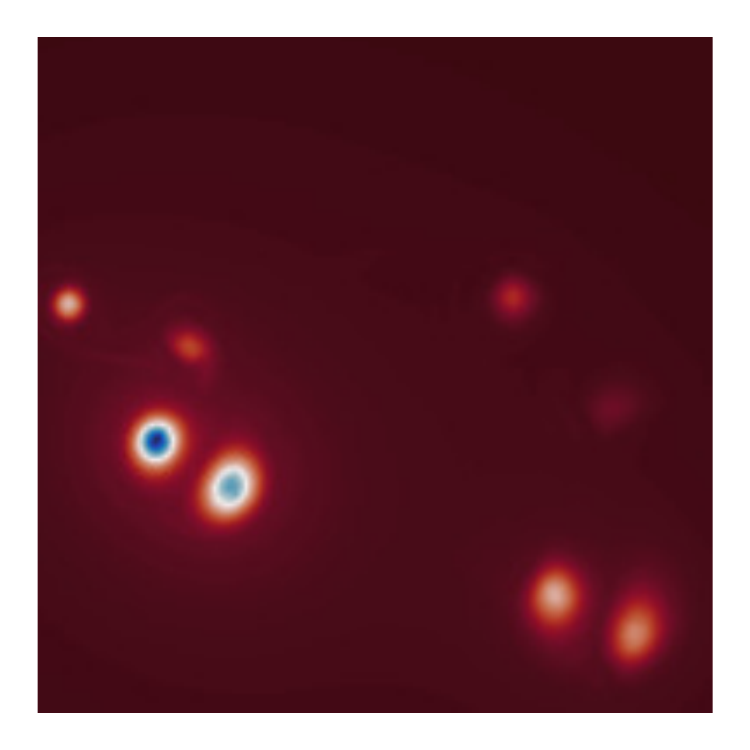}};
	\end{subfigure}		
	\\
	\vspace{-.06in}
	\hspace{.1in}	
	\begin{subfigure}[t]{0.315\textwidth}
		\centering
		\DeclareGraphicsExtensions{.pdf}
		 \hspace{-.1in}\tikz[baseline=(a.north)]\node[yscale=-1,inner sep=0,outer sep=0](a)
		{\begin{overpic}[width=1\textwidth]{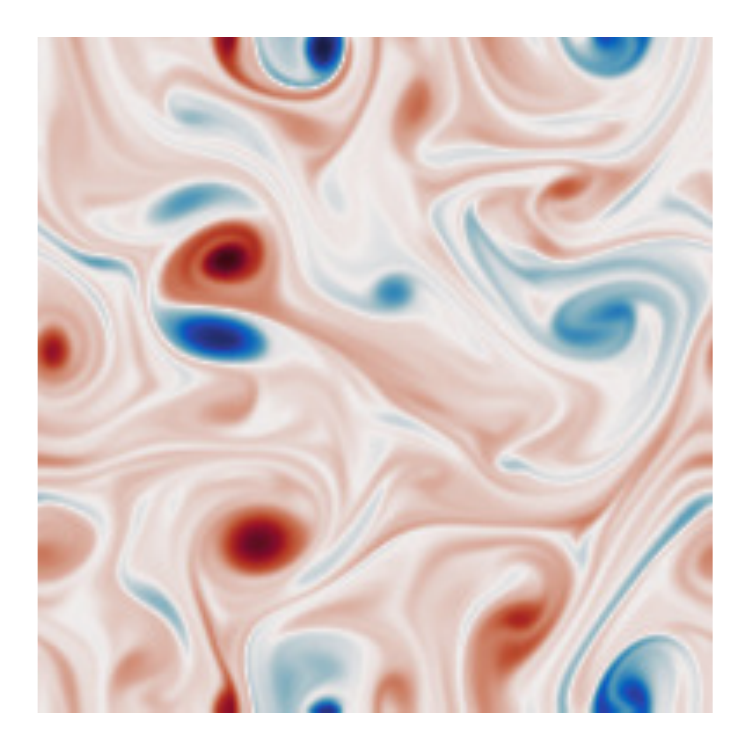}
		\put(-5,48){\reflectbox{\rotatebox[origin=c]{270}{2D Turbulence}}}
		\end{overpic}};
		\caption{Snapshots.}
	\end{subfigure}
	~
	\begin{subfigure}[t]{0.315\textwidth}
		\centering
		\DeclareGraphicsExtensions{.pdf}
		 \hspace{-.1in}\tikz[baseline=(a.north)]\node[yscale=-1,inner sep=0,outer sep=0](a)
		{\includegraphics[width=1\textwidth]{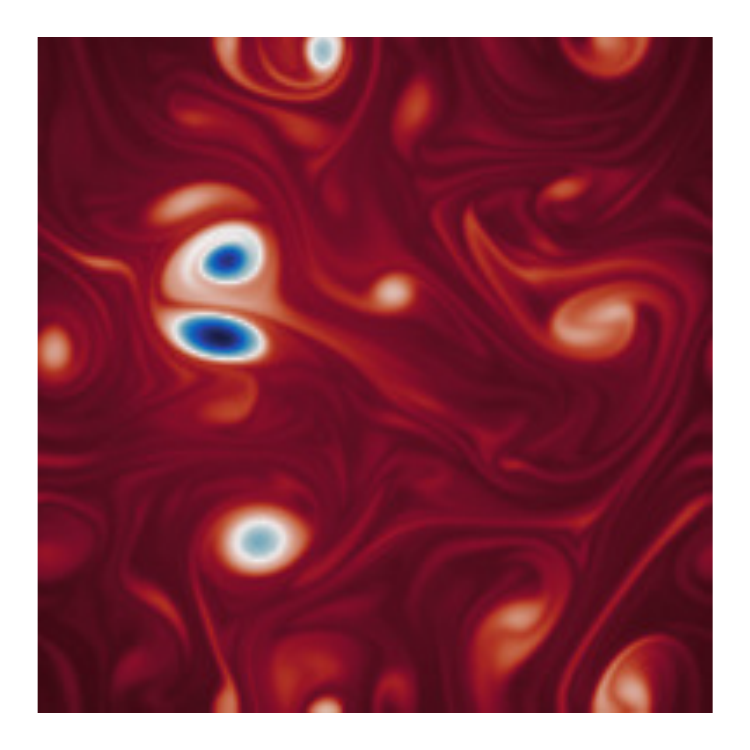}}; 
		\caption{Dominant eigenvectors.}
	\end{subfigure}
	~
	\begin{subfigure}[t]{0.315\textwidth}
		\centering
		\DeclareGraphicsExtensions{.pdf}
		 \hspace{-.1in}\tikz[baseline=(a.north)]\node[yscale=-1,inner sep=0,outer sep=0](a)
		{\includegraphics[width=1\textwidth]{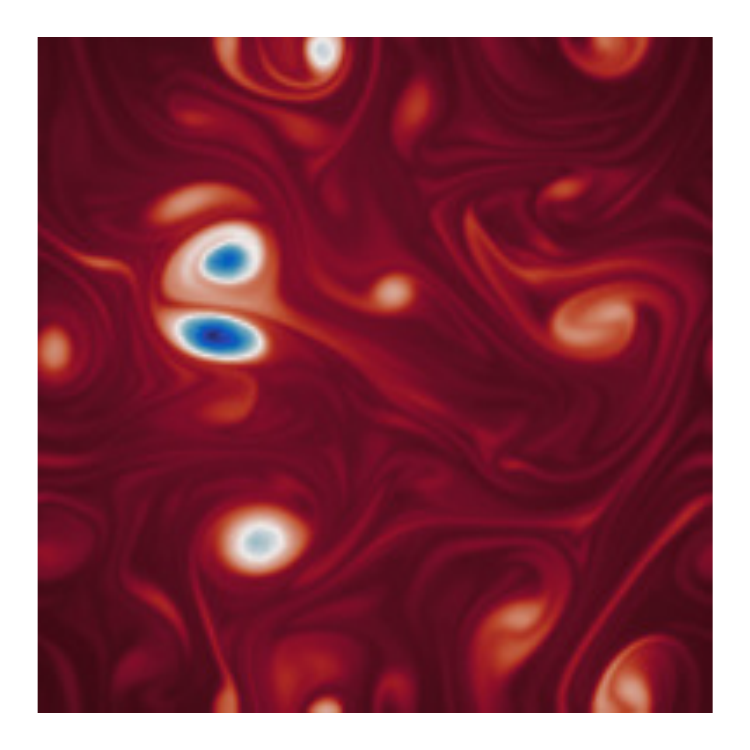}};
		\caption{Approximate eigenvectors.}
	\end{subfigure}	
	\vspace{.05in}	
	\caption{Qualitative comparison of the deterministic (b) and approximate (c) dominant eigenvector of the adjacency matrix generated from the flow field in (a). The top row shows the results for the airfoil wake  and the bottom shows the results for the two-dimensional isotropic turbulence. The approximation in (c) uses the Nystr\"om method with Halton sampling with 10\% of the columns of the adjacency matrix.}
	\label{fig:visual_results}
\end{figure}

\begin{figure}[!t]
\hspace{0.3in}
	\centering	
	\begin{subfigure}[t]{0.45\textwidth}
		\centering
		\DeclareGraphicsExtensions{.pdf}
		\hspace{-0.4in}\tikz[baseline=(a.north)]\node[yscale=-1,inner sep=0,outer sep=0](a)
		{\begin{overpic}[width=1\textwidth]{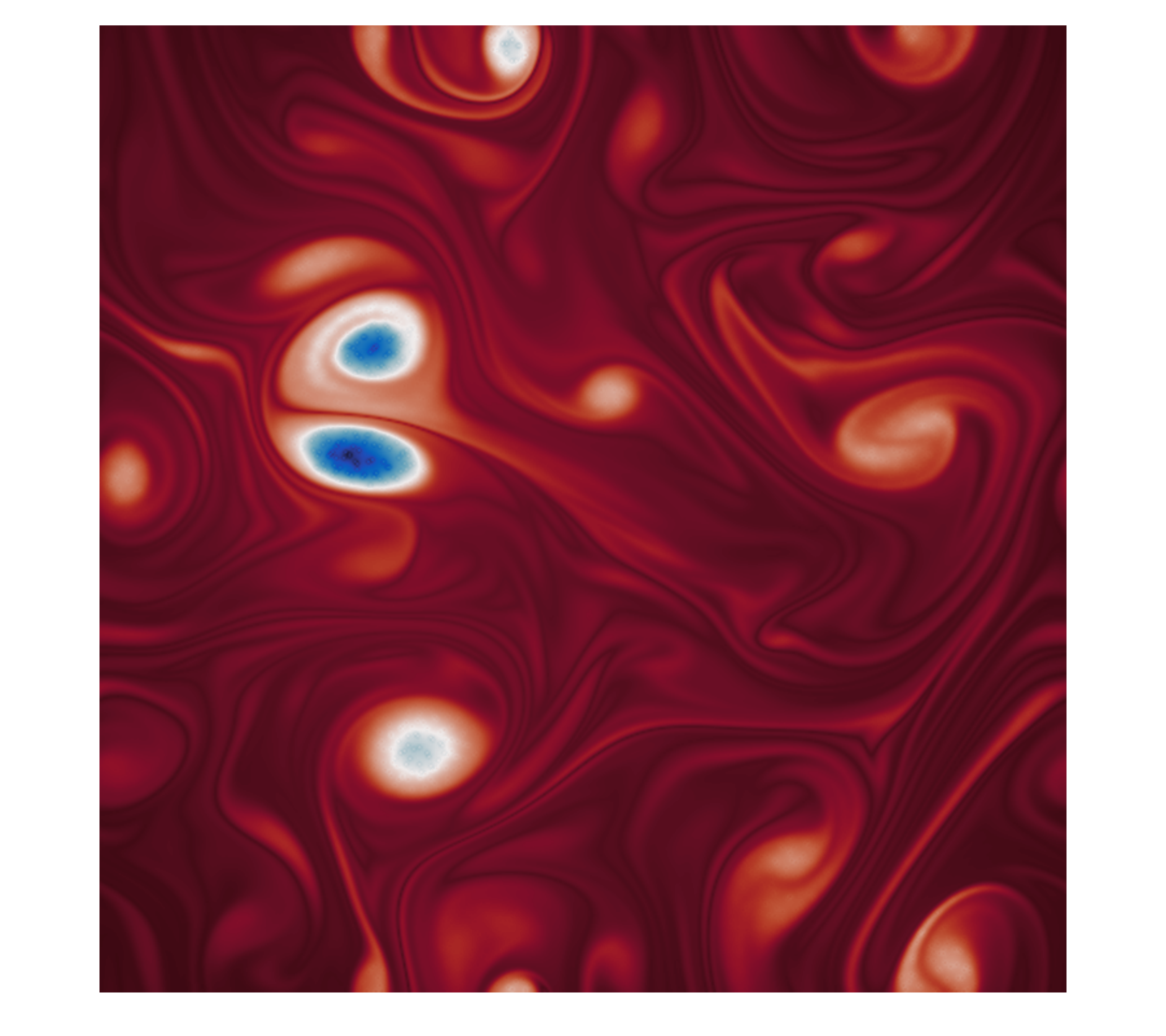}
		\end{overpic}};
		\caption{Eigenvector of dimension $512^2$.}
	\end{subfigure}
	~
	\begin{subfigure}[t]{0.45\textwidth}
		\centering
		\DeclareGraphicsExtensions{.pdf}
		\hspace{-0.4in}\tikz[baseline=(a.north)]\node[yscale=-1,inner sep=0,outer sep=0](a)
		{\begin{overpic}[width=1\textwidth]{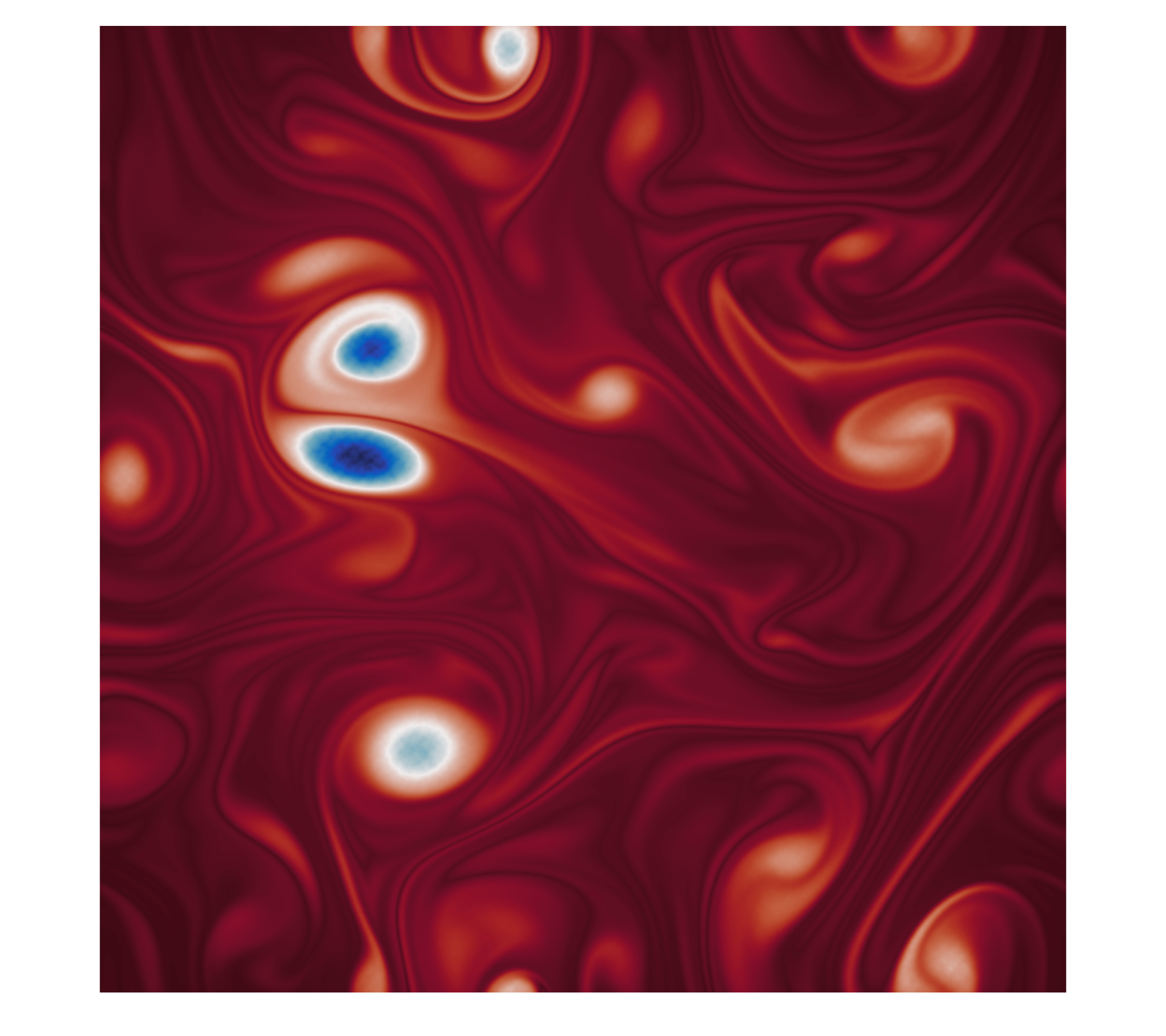}
		\end{overpic}};
		\caption{Eigenvector of dimension $1024^2$.}
	\end{subfigure}
	~
	\vspace{.05in}	
	\caption{Approximated leading eigenvectors of the adjacency matrices for higher-resolution isotropic flow fields. Here, we are using the Nystr\"om method with Halton sampling. By visual inspection we can see that sampling only about $10\%$ of the columns of the full adjacency matrix is sufficient for computing an approximated eigenvector that captures the dominant graph structure.}
	\label{fig:high_res}
\end{figure}

Figure~\ref{fig:visual_results} shows the leading eigenvectors of the adjacency matrix for the two flow examples, computed using deterministic and randomized methods.  
By visual inspection, the Nystr\"om method with quasi-random Halton sampling provides an excellent approximation of the leading eigenvector using only $10\%$ of the columns of the adjacency matrix. 
{ The resulting eigenvector centrality of the matrix $\mathbf{A}$ highlights strong vortex cores in both flows, as shown in Figure~\ref{fig:visual_results} (b). 
Vortices with lower strength have less centrality, emphasizing their decreased role in influencing the flow field.  
The 2P wake structure is clearly visible in the flow past an airfoil, and some shear layer structures are present in the 2D turbulent flow. 
The irrotational regions in both flows have the lowest centrality measure. 
These observations complement the traditional fluid dynamics  literature regarding highly influential structures in a vortical flow. The approximate eigenvector centrality measures capture both the highly influential vortex cores as well as the smaller, less influential  vortices, as seen in both flows. Moreover, the network-based measures take into account the spatial arrangements of the vortices to assess their influence.

In the computational comparison above, we only considered relatively low resolution flow fields because it was computationally prohibitive to compute the ground truth eigenvectors using deterministic methods.  
However, with randomized techniques it is possible to compute the leading eigenvector of the adjacency matrix for a $512\times 512$ resolution grid and a $1024\times 1024$ resolution grid, respectively, as shown in Figure~\ref{fig:high_res}. 
The adjacency matrix for the $1024\times 1024$ grid would require about $8$ TB of storage, if explicitly constructed. 
Instead, we sample only about $10\%$ of the columns to compute the Nystr\"om approximation. While we cannot quantify the approximation error, we can see by visual inspection that the eigenvector captures dominant flow structures.}

\begin{figure}[!t]
\hspace{0.14in}
 	\begin{subfigure}[t]{0.24\textwidth}
 	\centering
 		\hspace{-0.25in}\includegraphics[height=3.3cm]{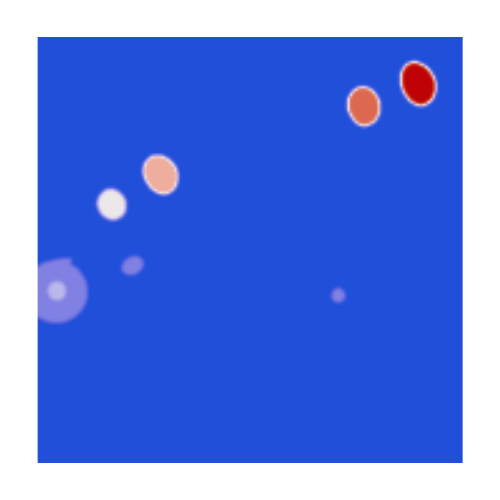}
 		\caption{Exact clusters.}
 	\end{subfigure}
 	\begin{subfigure}[t]{0.24\textwidth}
 	\centering
		\hspace{-0.3in}\includegraphics[height=3.3cm]{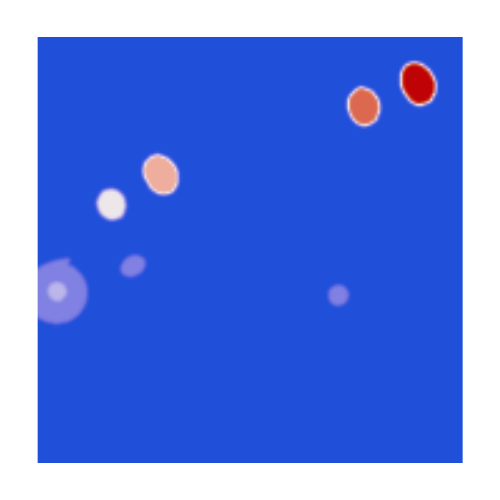}
		 \caption{Approximation.}
    \end{subfigure}
    \begin{subfigure}[t]{0.24\textwidth}
    \centering
 		\hspace{-0.3in}\includegraphics[height=3.2cm]{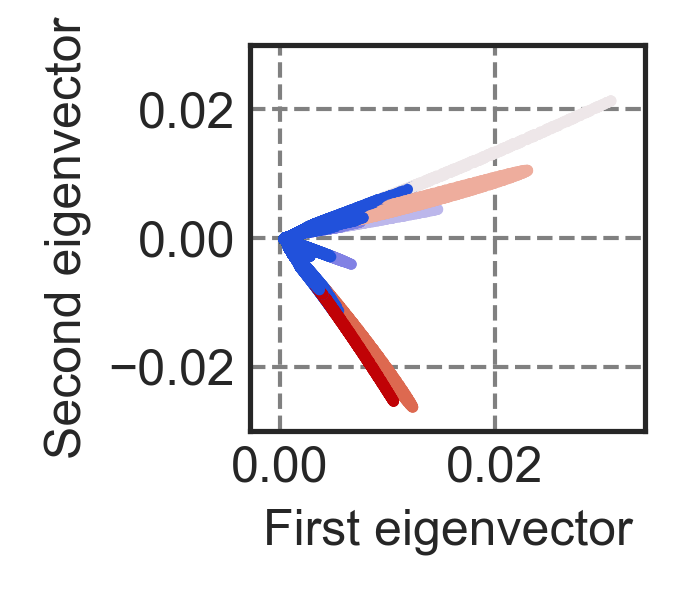}
 		 \caption{Scores (1 vs. 2).}
 	\end{subfigure}
 	\begin{subfigure}[t]{0.24\textwidth}
 	\centering
 		\hspace{-0.3in}\includegraphics[height=3.2cm]{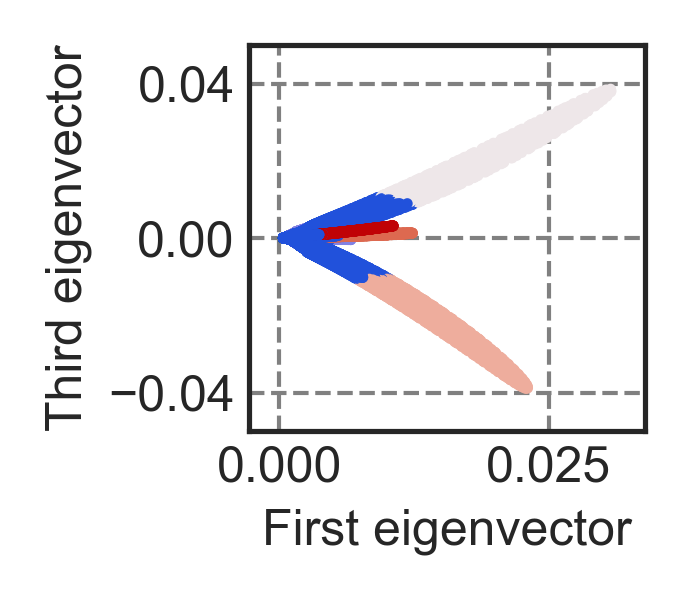}
 		 \caption{Scores (1 vs. 3).}
 	\end{subfigure}	
\vspace{.05in}	

	\caption{Spectral clustering using the top three eigenvectors of the adjacency matrix for the airfoil wake flow. In (a) we show the clusters that we found using the exact eigenvectors. It can be seen that the approximate eigenvectors allow us to reveal the same seven clusters, as shown in (b). In (c) and (d) we show the first three approximated eigenvector coordinates to visualize the clusters and subspaces. }
	\label{fig: clusters2}
\end{figure}

Finally, we investigate spectral clustering, which is one of the common uses of the leading eigenvectors of the adjacency matrix, providing a principled approach to cluster nodes into common communities.  
Figure~\ref{fig: clusters2} shows the results of spectral clustering based on the three leading eigenvectors obtained from the adjacency matrix.  
In both flow fields, the approximate clusters obtained using randomized methods closely resemble the clusters computed via deterministic eigendecomposition.   
The right panels in Figure~\ref{fig: clusters2} further shows the data plotted in these first three eigenvector coordinates to visualize the clusters and subspaces. 
{Typically, these leading eigenvectors are quite useful for determining relevant community structures within a network. 
	These communities have an important physical interpretation in vortical flow networks, hierarchically identifying and organizing vortex cores within the flow.  
	Because these communities are often determined from the few dominant eigenvectors, the randomized approach here can provide considerable computational advantages. 
	Further, since high-vorticity nodes have a large influence on the network interactions, it may be possible to improve convergence by preferential sampling using these community structures. 
}

\section{Discussion}
In this work, we have demonstrated the effective use of scalable algorithms from randomized linear algebra to accelerate network computations for large-scale fluid flows, as demonstrated on the wake behind an airfoil and two-dimensional isotropic turbulent flow.
In particular, we have used sampling-based methods to approximate the leading eigenvectors of the adjacency matrix of the vortical interaction network for cases where the data is so large that even single-pass algorithms are prohibitively expensive.  
Combining importance sampling, based on the probability distribution of detected communities, and the Nystr\"om method, the leading eigenvector can be accurately computed at a fraction of the computational cost and memory requirement, bypassing the need to construct or query the full $\mathbf{A}$ matrix.  
We also showed that quasi-random Halton sampling techniques outperform uniform sampling in both examples as they provide more comprehensive sampling coverage of the spatial domain. 

{There are a number of interesting future directions based on this work.  
Further study is required to apply these randomized network based analysis techniques to more complex three-dimensional turbulent flows. 
In addition, there are similar scaling issues in the related fields of almost invariant sets and set-oriented methods~\cite{Dellnitz2001book,Dellnitz2002hds,Froyland2009pd,Froyland2010chaos,Tallapragada2013}, where the state-transition operator may be viewed as a graph on the high-dimensional state-space. 
In the transfer operator case, machine learning has been used to identify a data-driven discretization of the state-space into clusters, which serve as nodes that scale with the intrinsic rank of the attractor, rather than the high-dimensional measurement dimension~\cite{Kaiser2014jfm}; similar clustering has been used to identify subspaces for POD--Galerkin models~\cite{Amsallem2012ijnme}. 
These methods have been used to determine eigenvectors of the Perron--Frobenius operator, the dual of the Koopman operator, establishing a strong connection between sensitivity and coherence in fluid flows.  
It will be interesting to mathematically connect those approaches with the randomized methods considered in the present study.  }

\section*{Acknowledgments}
We acknowledge support from the Army Research Office (W911NF-17-1-0118) and the Air Force Office of Scientific Research (FA9550-16-1-0650). ZB also acknowledges support by the U.S. Department of Energy, Office of Science, under Award Number DE-AC02-05CH11231. NBE would also like to acknowledge Amazon Web Services for supporting the
project with EC2 credits. 
We would like to thank Bing Brunton, Nathan Kutz, Jean-Christophe Loiseau, Krithika Manohar, Aditya Nair, and Bernd Noack for valuable discussions. 

\begin{spacing}{.9}
\small
\setlength{\bibsep}{6.5pt}
\bibliographystyle{ieeetr}

\end{spacing}
\end{document}